\documentclass[thmsa]{article}

\usepackage{amsmath,amsthm,amssymb}
\usepackage{amsfonts}
\usepackage{graphicx}

\newtheorem{theorem}{Theorem}
\newtheorem{corollary}[theorem]{Corollary}
\newtheorem{lemma}[theorem]{Lemma}

\newtheorem{example}[theorem]{Example}
\newtheorem{remark}[theorem]{Remark}

\def\proof{\noindent{{\bf Proof. }}}
\def\endproof{ $\blacksquare$}

\makeatletter
\def\@makefnmark{}
\makeatother

\DeclareMathOperator*{\sgn}{sgn}

\begin{document}

\title{Nonlinear Schr\"{o}dinger equations without compatibility \\ conditions on the potentials}

\author{
Michela Guida\thanks{%
Dipartimento di Matematica ``Giuseppe Peano'', Universit\`{a} degli Studi di
Torino, Via Carlo Alberto 10, 10123 Torino, Italy. e-mail: \texttt{%
michela.guida@unito.it}},\quad 
Sergio Rolando\thanks{%
Dipartimento di Matematica e Applicazioni, Universit\`{a} di Milano-Bicocca,
Via Roberto Cozzi 53, 20125 Milano, Italy. e-mail: \texttt{%
sergio.rolando@unito.it}}}
\date{\vspace{-0.7cm}}
\maketitle

\begin{abstract}
We study the existence of nonnegative solutions (and ground states) to
nonlinear Schr\"{o}dinger equations in $\mathbb{R}^{N}$ with radial potentials
and super-linear or sub-linear nonlinearities. The potentials satisfy power
type estimates at the origin and at infinity, but no compatibility condition
is required on their growth (or decay) rates at zero and infinity. In this
respect our results extend some well known results in the literature and we
also believe that they can highlight the role of the sum of Lebesgue spaces
in studying nonlinear equations with weights.\bigskip

\noindent \textbf{MSC (2010):} Primary 35J60; Secondary 35J05, 35J20,
46E30\smallskip

\noindent \textbf{Keywords:} Nonlinear Schr\"{o}dinger equation, unbounded
or decaying potentials, sum of weighted Lebesgue speces, ground
states\medskip
\end{abstract}

\footnotetext{
The authors are members of the Gruppo Nazionale di Alta Matematica (INdAM).}%

\section{Introduction and main results}

We study the semilinear elliptic equation 
\begin{equation}
-\triangle u+V\left( \left| x\right| \right) u=K\left( \left| x\right|
\right) f\left( u\right) \quad \textrm{in }\mathbb{R}^{N},~N\geq 3,  \label{Eq.}
\end{equation}
where $f:\mathbb{R}\rightarrow \mathbb{R}$ is a continuous function such that $%
f\left( 0\right) =0$ and $V$, $K$ satisfy the following assumptions:

\begin{itemize}
\item[$\left( \mathbf{V}\right) $]  $V:\left( 0,+\infty \right) \rightarrow
\left[ 0,+\infty \right) $ is a continuous function such that 
\[
\liminf_{r\rightarrow 0^{+}}\frac{V\left( r\right) }{r^{a_{0}}}>0\quad \textrm{%
and}\quad \liminf_{r\rightarrow +\infty }\frac{V\left( r\right) }{r^{a}}%
>0\quad \textrm{for some }a_{0},a\in \mathbb{R};
\]

\item[$\left( \mathbf{K}\right) $]  $K:\left( 0,+\infty \right) \rightarrow
\left( 0,+\infty \right) $ is a continuous function such that 
\[
\limsup_{r\rightarrow 0^{+}}\frac{K\left( r\right) }{r^{b_{0}}}<\infty \quad 
\textrm{and}\quad \limsup_{r\rightarrow +\infty }\frac{K\left( r\right) }{r^{b}%
}<\infty \quad \textrm{for some }b_{0},b\in \mathbb{R}.
\]
\newpage 
\end{itemize}

\noindent More precisely, we are interested in finding nontrivial
nonnegative radial solutions in the following weak sense (see also Remark 
\ref{RMK: finale}.\ref{RMK:symm-crit}): we call \emph{radial solution} to
Eq. (\ref{Eq.}) any $u\in H_{V,\mathrm{r}}^{1}$ such that 
\begin{equation}
\int_{\mathbb{R}^{N}}\nabla u\cdot \nabla h\,dx+\int_{\mathbb{R}^{N}}V\left(
\left| x\right| \right) uh\,dx=\int_{\mathbb{R}^{N}}K\left( \left| x\right|
\right) f\left( u\right) h\,dx\qquad \textrm{for all }h\in H_{V,\mathrm{r}%
}^{1},  \label{weak solution}
\end{equation}
where 
\begin{equation}
H_{V,\mathrm{r}}^{1}=H_{V,\mathrm{r}}^{1}\left( \mathbb{R}^{N}\right) :=\left\{
u\in H_{V}^{1}\left( \mathbb{R}^{N}\right) :u\left( x\right) =u\left( \left|
x\right| \right) \right\}  \label{H^1_V,r}
\end{equation}
is the radial subspace of 
\begin{equation}
H_{V}^{1}=H_{V}^{1}\left( \mathbb{R}^{N}\right) :=\left\{ u\in D^{1,2}\left( 
\mathbb{R}^{N}\right) :\int_{\mathbb{R}^{N}}V\left( \left| x\right| \right)
u^{2}dx<\infty \right\} .  \label{H^1_V}
\end{equation}
Here $D^{1,2}(\mathbb{R}^{N})=\{u\in L^{2^{*}}(\mathbb{R}^{N}):\left| \nabla
u\right| \in L^{2}(\mathbb{R}^{N})\}$, $2^{*}:=2N/(N-2)$, denotes the usual
Sobolev space, which identifies with the completion of $C_{\mathrm{c}%
}^{\infty }(\mathbb{R}^{N})$ with respect to the $L^{2}$ norm of the
gradient. Of course, $u\left( x\right) =u\left( \left| x\right| \right) $
means that $u$ is invariant under the action on $H_{V}^{1}$ of the
orthogonal group of $\mathbb{R}^{N}$.

By well known arguments, the nonnegative weak solutions to Eq. (\ref{Eq.})
lead to special solutions (\textit{solitary waves} and \textit{solitons})
for several nonlinear field theories, such as nonlinear Schr\"{o}dinger and
Klein-Gordon equations, which arise in many branches of mathematical
physics, such as nonlinear optics, plasma physics, condensed matter physics
and cosmology (see e.g. \cite{BBR1,BFmonograph,YangY}). In this respect,
since the early studies of \cite{Beres-Lions,Floer-Wein,Rabi92,Strauss}, Eq.
(\ref{Eq.}) has been massively addressed in the mathematical literature,
recently focusing on the case of $V$ possibly vanishing at infinity, that
is, $\liminf_{\left| x\right| \rightarrow \infty }V\left( \left| x\right|
\right) =0$ (some first results on such a case can be found in \cite
{Ambr-Fel-Malch,BR,Be-Gr-Mic,Be-Gr-Mic 2}; for more recent bibliography, see
e.g. \cite{Alves-Souto-13,BGRnonex,BR TMA,Catrina nonex,SuTian12} and the
references therein).

The most recent and general existence results for radial solutions to Eq. (%
\ref{Eq.}) under assumptions $\left( \mathbf{V}\right) $ and $\left( \mathbf{%
K}\right) $, unifying and extending the previously existing ones, are
contained in \cite{Su-Wang-Will p} and \cite{SuTian12}, which respectively
concern the case of super-linear and sub-linear nonlinearities.

The result of \cite{Su-Wang-Will p}, rewritten in a suitable form for
comparing with our results, is Theorem \ref{THM: SWW} below, which uses the
following notation. For every $a_{0}\in \mathbb{R}$, set 
\begin{equation}
\underline{b}\left( a_{0}\right) :=\left\{ 
\begin{array}{ll}
-\infty & \textrm{if }a_{0}<-\left( 2N-2\right) \\ 
\min \left\{ a_{0},-2\right\} ~ & \textrm{if }a_{0}\geq -\left( 2N-2\right) .
\end{array}
\right.  \label{b_:=}
\end{equation}
Then, for $a,b,a_{0}\in \mathbb{R}$ and $b_{0}>\underline{b}\left( a_{0}\right) 
$, define the functions 
\[
\underline{q}=\underline{q}\left( a,b,a_{0},b_{0}\right) :=\left\{ 
\begin{array}{lll}
\max \left\{ 2,\,2\frac{N+b}{N-2}\right\} & ~\smallskip & \textrm{if }a\leq
-2,~b_{0}>\min \left\{ -2,a_{0}\right\} \\ 
\max \left\{ 2,\,2\frac{N+b}{N-2},\,2\frac{2N-2+2b_{0}-a_{0}}{2N-2+a_{0}}%
\right\} & \smallskip & \textrm{if }a\leq -2,~b_{0}\leq a_{0}<-\left(
2N-2\right) \\ 
\max \left\{ 2,\,2\frac{2N-2+2b-a}{2N-2+a}\right\} & \smallskip & \textrm{if }%
a>-2,~b_{0}>\min \left\{ -2,a_{0}\right\} \\ 
\max \left\{ 2,\,2\frac{2N-2+2b-a}{2N-2+a},\,2\frac{2N-2+2b_{0}-a_{0}}{%
2N-2+a_{0}}\right\} &  & \textrm{if }a>-2,~b_{0}\leq a_{0}<-\left( 2N-2\right)
\end{array}
\right. 
\]
and 
\[
\overline{q}=\overline{q}\left( a_{0},b_{0}\right) :=\left\{ 
\begin{array}{lll}
+\infty & ~\smallskip & \textrm{if }a_{0}<-\left( 2N-2\right) ~\textrm{or}%
~~a_{0}=-\left( 2N-2\right) <b_{0} \\ 
2\frac{2N-2+2b_{0}-a_{0}}{2N-2+a_{0}} & \smallskip & \textrm{if }-\left(
2N-2\right) <a_{0}<-2,~b_{0}>a_{0} \\ 
2\frac{N+b_{0}}{N-2} &  & \textrm{if }a_{0}\geq -2,~b_{0}>-2.
\end{array}
\right. 
\]
Observe that one always has $\underline{q}\geq 2$ and $\overline{q}>2$.

\begin{theorem}[{\cite[Theorem 5]{Su-Wang-Will p}}]
\label{THM: SWW}Assume $\left( \mathbf{V}\right) ,\left( \mathbf{K}\right) $
with $a_{0},a,b\in \mathbb{R}$ and $b_{0}>\underline{b}\left( a_{0}\right) $.
Assume furthermore that $\underline{q}<\overline{q}$. Then Eq. (\ref{Eq.})
has a nonnegative nontrivial radial solution for every continuous $f:\mathbb{R}%
\rightarrow \mathbb{R}$ satisfying:

\begin{itemize}
\item[$\left( \mathbf{f}_{1}\right) $]  
$\displaystyle \sup_{t>0}\,\frac{\left| f\left( t\right) \right| }{t^{q-1}}<+\infty $ for
some $q\in (\underline{q},\overline{q});$

\item[$\left( \mathbf{f}_{2}\right) $]  
$\exists \theta >2$ such that $0<\theta F\left( t\right) \leq f\left( t\right) t$ for all $t\in \mathbb{R}.$
\end{itemize}
\end{theorem}

In $\left( \mathbf{f}_{2}\right) $ and everywhere in the following, we
denote $F\left( t\right) :=\int_{0}^{t}f\left( s\right) ds$.

\begin{remark}
To be precise, instead of $\left( \mathbf{f}_{1}\right) $, the growth
condition used in \cite[Theorem 5]{Su-Wang-Will p} is 
\begin{equation}
\sup_{t\in \mathbb{R}}\,\frac{\left| f\left( t\right) \right| }{\left| t\right|
^{q_{1}-1}+\left| t\right| ^{q_{2}-1}}<+\infty \quad \textrm{for some }%
q_{1},q_{2}\in \left( \underline{q},\overline{q}\right) ,  \label{growth sum}
\end{equation}
but the difference between (\ref{growth sum}) and $\left( \mathbf{f}%
_{1}\right) $ is not essential. Indeed, we can just let $t>0$ in (\ref
{growth sum}) because we deal with nonnegative solutions, and the use of a
sum of powers is a standard generalization of $\left( \mathbf{f}_{1}\right) $.
\end{remark}

In order to recall the existence result of \cite{SuTian12}, we need some
further notation. Define the following subsets of $\mathbb{R}^{2}$: 
\[
\begin{array}{ll}
\mathcal{A}_{1}:=\left\{ \left( a,b\right) :\max \left\{ -\frac{N+2}{2},%
\frac{a-2}{2}\right\} \leq b<-2\right\} ,\medskip & \mathcal{B}_{1}:=\left\{
\left( a_{0},b_{0}\right) :\max \left\{ -\frac{N+2}{2},\frac{a_{0}-2}{2}%
\right\} <b_{0}\leq -2\right\} , \\ 
\mathcal{A}_{2}:=\left\{ \left( a,b\right) :-\frac{N+2}{2}\leq b<\min
\left\{ -2,\frac{a-2N-2}{4}\right\} \right\} ,\medskip & \mathcal{B}%
_{2}:=\left\{ \left( a_{0},b_{0}\right) :-\frac{N+2}{2}<b_{0}\leq -2\leq
a_{0}\right\} , \\ 
\mathcal{A}_{3}:=\left\{ \left( a,b\right) :a\leq -2,~-\frac{N+2}{2}<b<\frac{%
a-2}{2}\right\} ,\medskip & \mathcal{B}_{3}:=\left\{ \left(
a_{0},b_{0}\right) :a_{0}<-2,~-\frac{N+2}{2}<b_{0}\leq \frac{a_{0}-2}{2}%
\right\} , \\ 
\mathcal{A}_{4}:=\left\{ \left( a,b\right) :b\leq -\frac{N+2}{2},~\frac{%
a-2N-2}{4}\leq b<\frac{a-2}{2}\right\} ,\medskip & \mathcal{B}_{4}:=\left\{
\left( a_{0},b_{0}\right) :b_{0}<-\frac{N+2}{2},~\frac{a_{0}-2N-2}{4}%
<b_{0}\leq \frac{a_{0}-2}{2}\right\} , \\ 
\mathcal{A}_{5}:=\left\{ \left( a,b\right) :a>-2,~\frac{a-2N-2}{4}\leq b<%
\frac{a-2}{2}\right\} ,\medskip & \mathcal{B}_{5}:=\left\{ \left(
a_{0},b_{0}\right) :b_{0}\geq -2,~\frac{a_{0}-2N-2}{4}<b_{0}\leq \frac{%
a_{0}-2}{2}\right\} , \\ 
\mathcal{B}:=\mathcal{B}_{1}\cup ...\cup \mathcal{B}_{5}, & \mathcal{B}%
_{6}:=\left\{ \left( a_{0},b_{0}\right) :\frac{a_{0}-2}{2}<b_{0}\leq \frac{%
a_{0}-2N-2}{4}\right\} .
\end{array}
\]
Then, for $\left( a,b\right) \in \mathcal{A}_{1}\cup ...\cup \mathcal{A}_{5}$
and $\left( a_{0},b_{0}\right) \in \mathcal{B}\cup \mathcal{B}_{6}$, define
the functions 
\[
\underline{\underline{q}}=\underline{\underline{q}}\left(
a,b,a_{0},b_{0}\right) :=\left\{ 
\begin{array}{lll}
2\frac{N+b}{N-2} & ~\smallskip & \textrm{if }\left( a,b\right) \in \mathcal{A}%
_{1}\cup \mathcal{A}_{2}\cup \mathcal{A}_{3},~\left( a_{0},b_{0}\right) \in 
\mathcal{B} \\ 
\max \left\{ 2\frac{N+b}{N-2},\,4\frac{N+b_{0}}{2N-2+a_{0}}\right\} & 
\smallskip & \textrm{if }\left( a,b\right) \in \mathcal{A}_{1}\cup \mathcal{A}%
_{2}\cup \mathcal{A}_{3},~\left( a_{0},b_{0}\right) \in \mathcal{B}_{6} \\ 
4\frac{N+b}{2N-2+a} & \smallskip & \textrm{if }\left( a,b\right) \in \mathcal{A%
}_{4}\cup \mathcal{A}_{5},~\left( a_{0},b_{0}\right) \in \mathcal{B} \\ 
\max \left\{ 4\frac{N+b}{2N-2+a},\,4\frac{N+b_{0}}{2N-2+a_{0}}\right\} &  & 
\textrm{if }\left( a,b\right) \in \mathcal{A}_{4}\cup \mathcal{A}_{5},~\left(
a_{0},b_{0}\right) \in \mathcal{B}_{6}
\end{array}
\right. 
\]
and 
\[
\overline{\overline{q}}=\overline{\overline{q}}\left( a_{0},b_{0}\right)
:=\left\{ 
\begin{array}{lll}
2\frac{N+b_{0}}{N-2} & ~\smallskip & \textrm{if }\left( a_{0},b_{0}\right) \in 
\mathcal{B}_{1}\cup \mathcal{B}_{2} \\ 
4\frac{N+b_{0}}{2N-2+a_{0}} & \smallskip & \textrm{if }\left(
a_{0},b_{0}\right) \in \mathcal{B}_{3}\cup \mathcal{B}_{4}\cup \mathcal{B}%
_{5} \\ 
2 &  & \textrm{if }\left( a_{0},b_{0}\right) \in \mathcal{B}_{6}.
\end{array}
\right. 
\]
Observe that one always has $1\leq \underline{\underline{q}}<2$ and $1<%
\overline{\overline{q}}\leq 2$.

\begin{theorem}[{\cite[Theorem 1.3]{SuTian12}}]
\label{THM: ST}Assume $\left( \mathbf{V}\right) ,\left( \mathbf{K}\right) $
with $\left( a,b\right) \in \mathcal{A}_{1}\cup ...\cup \mathcal{A}_{5}$ and 
$\left( a_{0},b_{0}\right) \in \mathcal{B}\cup \mathcal{B}_{6}$. Assume
furthermore that $\underline{\underline{q}}<\overline{\overline{q}}$. Then
Eq. (\ref{Eq.}) with $f\left( u\right) =\left| u\right| ^{q-2}u$ has a
nonnegative nontrivial radial solution provided that $q\in (\underline{%
\underline{q}},\overline{\overline{q}})$.
\end{theorem}

The spirit of the above Theorems \ref{THM: SWW} and \ref{THM: ST} is
essentially the following: a compatibility condition between the behaviours
of the potentials at zero and at infinity is required ($\underline{q}<%
\overline{q}$ or $\underline{\underline{q}}<\overline{\overline{q}}$) and a
solution is then provided if the nonlinearity grows compatibly with the
potentials ($q$ between $\underline{q}$, $\overline{q}$ or $\underline{%
\underline{q}}$, $\overline{\overline{q}}$).

Here we still require some compatibility between the nonlinearity and the
potentials, but we remove any compatibility assumption between how the
potentials behave at zero and infinity, getting existence results that
contain and extend Theorems \ref{THM: SWW} and \ref{THM: ST}.

In order to state our results, we need some preliminary notations, which
essentially consist in defining two intervals $\mathcal{I}_{1}$ and $%
\mathcal{I}_{2}$ which will provide a way of expressing the compatibility
between the nonlinearity and the potentials required in order to get
existence. In this respect, such intervals play the same role of the
limiting exponents $\underline{q}$, $\overline{q}$, $\underline{\underline{q}%
}$, $\overline{\overline{q}}$ of Theorems \ref{THM: SWW} and \ref{THM: ST},
but in an unified way for both the cases of super-linear and sub-linear
nonlinearities (see assumptions $\left( \mathbf{f}_{3}\right) $ and $\left( 
\mathbf{f}_{7}\right) $ below). On a more technical level, $\mathcal{I}_{1}$
and $\mathcal{I}_{2}$ are the exact ranges of exponents for which we can
prove the compactness result given in Lemma \ref{LEM: comp embedd} of
Section \ref{SEC: proofs}.

For every $a_{0}\in \mathbb{R}$, define 
\[
b_{*}\left( a_{0}\right) :=\left\{ 
\begin{array}{lll}
-\infty & ~\smallskip & \textrm{if }a_{0}<-\left( 2N-2\right) \\ 
\min \left\{ a_{0},-\frac{N-a_{0}}{2},-\frac{N+2}{2}\right\} &  & \textrm{if }%
a_{0}\geq -\left( 2N-2\right) .
\end{array}
\right. 
\]
Recalling definition (\ref{b_:=}) of $\underline{b}$, observe that $%
b_{*}\left( a_{0}\right) \leq \underline{b}\left( a_{0}\right) $ for every $%
a_{0}\in \mathbb{R}$ (precisely: $b_{*}=\underline{b}$ for $a_{0}\leq -N$ and $%
b_{*}<\underline{b}$ for $a_{0}>-N$). Then, for $a_{0}\in \mathbb{R}$ and $%
b_{0}>b_{*}\left( a_{0}\right) $, define the functions 
\[
q_{*}\left( a_{0},b_{0}\right) :=\left\{ 
\begin{array}{lll}
\max \left\{ 1,\,2\frac{N+b_{0}}{N+a_{0}},\,2\frac{2N-2+2b_{0}-a_{0}}{%
2N-2+a_{0}}\right\} & ~\smallskip & \textrm{if }a_{0}<-\left( 2N-2\right) \\ 
\max \left\{ 1,\,2\frac{N+b_{0}}{N+a_{0}}\right\} & \smallskip & \textrm{if }%
-\left( 2N-2\right) \leq a_{0}<-N \\ 
1 &  & \textrm{if }a_{0}\geq -N,
\end{array}
\right. 
\]
\[
q^{*}\left( a_{0},b_{0}\right) :=\left\{ 
\begin{array}{lll}
+\infty & ~\smallskip & \textrm{if }a_{0}\leq -\left( 2N-2\right) \\ 
2\frac{2N-2+2b_{0}-a_{0}}{2N-2+a_{0}} & \smallskip & \textrm{if }-\left(
2N-2\right) <a_{0}\leq -N \\ 
\min \left\{ 2\frac{N+b_{0}}{N+a_{0}},\,2\frac{2N-2+2b_{0}-a_{0}}{2N-2+a_{0}}%
\right\} & \smallskip & \textrm{if }-N<a_{0}<-2 \\ 
2\frac{N+b_{0}}{N-2} &  & \textrm{if }a_{0}\geq -2
\end{array}
\right. 
\]
and the interval 
\begin{equation}
\mathcal{I}_{1}=\mathcal{I}_{1}\left( a_{0},b_{0}\right) :=\left(
q_{*}\left( a_{0},b_{0}\right) ,\,q^{*}\left( a_{0},b_{0}\right) \right) .
\label{I_1:=}
\end{equation}
Note that $b_{0}>b_{*}\left( a_{0}\right) $ is equivalent to $q_{*}\left(
a_{0},b_{0}\right) <q^{*}\left( a_{0},b_{0}\right) $, i.e., $\mathcal{I}%
_{1}\neq \varnothing $. Finally, for every $a,b\in \mathbb{R}$, define the
function 
\[
q_{**}\left( a,b\right) :=\left\{ 
\begin{array}{lll}
\max \left\{ 1,\,2\frac{N+b}{N-2}\right\} & ~\smallskip & \textrm{if }a\leq -2
\\ 
\max \left\{ 1,\,2\frac{N+b}{N+a},\,2\frac{2N-2+2b-a}{2N-2+a}\right\} &  & 
\textrm{if }a>-2
\end{array}
\right. 
\]
and the interval 
\begin{equation}
\mathcal{I}_{2}=\mathcal{I}_{2}\left( a,b\right) :=\left( q_{**}\left(
a,b\right) ,+\infty \right) .  \label{I_2:=}
\end{equation}
In order to ease the visualization of the intervals $\mathcal{I}_{1}$ and $%
\mathcal{I}_{2}$, the graphs of the functions $q_{*}\left( a_{0},\cdot
\right) $, $q^{*}\left( a_{0},\cdot \right) $ and $q_{**}\left( a,\cdot
\right) $, with $a_{0}$ and $a$ fixed to different meaningful values, are
plotted in Figures 1-8 below.
\bigskip\bigskip

$
\begin{tabular}[t]{l}
\begin{tabular}{l}
\includegraphics[height=1.50in]{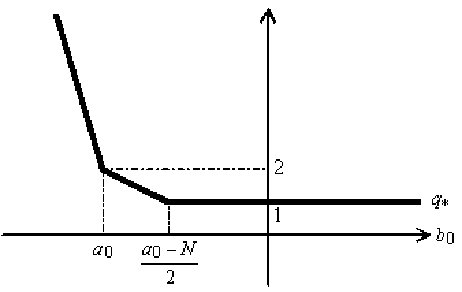}
\smallskip  \\ 
\textbf{Fig.1.} $q_{*}\left( a_{0},\cdot \right) $ for $a_{0}<-\left(
2N-2\right) \smallskip $ \\ 
($q^{*}\left( a_{0},\cdot \right) =+\infty $)
\end{tabular}
\end{tabular}
\begin{tabular}[t]{l}
\begin{tabular}{l}
\includegraphics[height=1.50in]{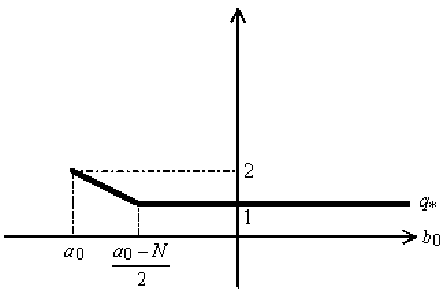}
\smallskip  \\ 
\textbf{Fig.2.} $q_{*}\left( a_{0},\cdot \right) $ for $a_{0}=-\left(
2N-2\right) \smallskip $ \\ 
($q^{*}\left( a_{0},\cdot \right) =+\infty $)
\end{tabular}
\end{tabular}
$\bigskip\bigskip

$
\begin{tabular}[t]{l}
\begin{tabular}{l}
\includegraphics[height=1.50in]{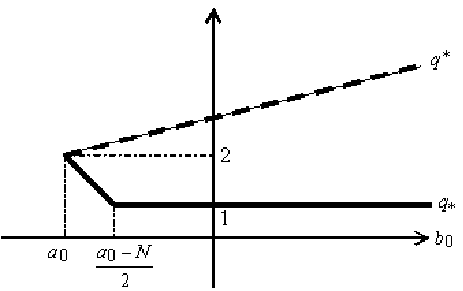}
\smallskip  \\ 
\textbf{Fig.3.} $q_{*}\left( a_{0},\cdot \right) $ and $q^{*}\left(
a_{0},\cdot \right) $ $\smallskip $ \\ 
for $-\left( 2N-2\right) <a_{0}<-N$%
\end{tabular}
\end{tabular}
~ 
\begin{tabular}[t]{l}
\begin{tabular}{l}
\includegraphics[height=1.50in]{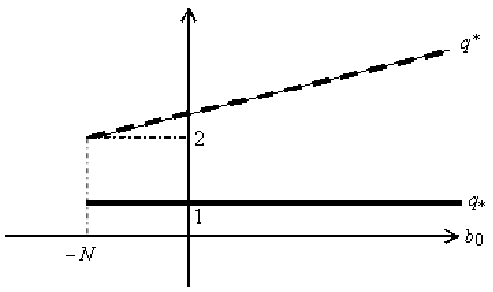}
\smallskip  \\ 
\textbf{Fig.4.} $q_{*}\left( a_{0},\cdot \right) $ and $q^{*}\left(
a_{0},\cdot \right) $ $\smallskip $ \\ 
for $a_{0}=-N$%
\end{tabular}
\end{tabular}
$\bigskip\bigskip

$
\begin{tabular}[t]{l}
\begin{tabular}{l}
\includegraphics[height=1.50in]{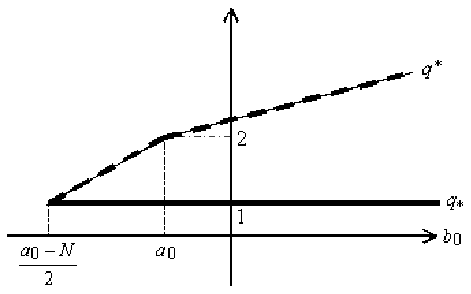}
\smallskip \\ 
\textbf{Fig.5.} $q_{*}\left( a_{0},\cdot \right) $ and $q^{*}\left(
a_{0},\cdot \right) $ $\smallskip $ \\ 
for $-N<a_{0}<-2$%
\end{tabular}
\end{tabular}
\begin{tabular}[t]{l}
\begin{tabular}{l}
\includegraphics[height=1.50in]{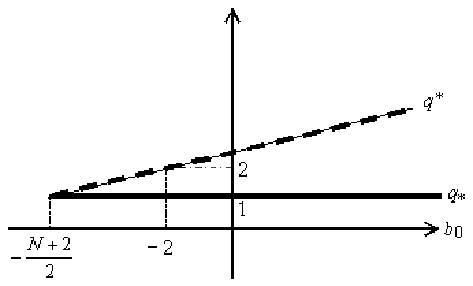}
\smallskip  \\ 
\textbf{Fig.6.} $q_{*}\left( a_{0},\cdot \right) $ and $q^{*}\left(
a_{0},\cdot \right) $ $\smallskip $ \\ 
for $a_{0}\geq -2$%
\end{tabular}
\end{tabular}
$
\pagebreak

$
\begin{tabular}[t]{l}
\begin{tabular}{l}
\includegraphics[height=1.50in]{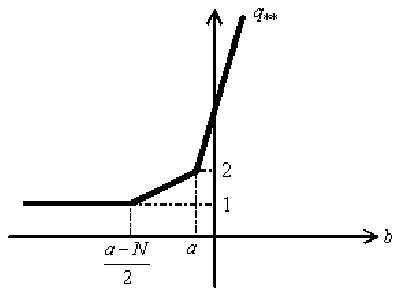}
\smallskip  \\ 
\textbf{Fig.7.} $q_{**}\left( a,\cdot \right) $ for $a\leq -2$%
\end{tabular}
\end{tabular}
~\,\qquad 
\begin{tabular}[t]{l}
\begin{tabular}{l}
\includegraphics[height=1.50in]{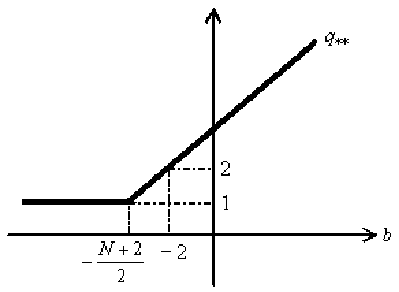}
\smallskip  \\ 
\textbf{Fig.8.} $q_{**}\left( a,\cdot \right) $ for $a>-2$%
\end{tabular}
\end{tabular}
$\bigskip\bigskip

Notice that:

\begin{itemize}
\item  both $\mathcal{I}_{1}$ and $\mathcal{I}_{2}$ are contained in $\left(
1,+\infty \right) $ for every $a_{0},a,b\in \mathbb{R}$ and $b_{0}>b_{*}\left(
a_{0}\right) $;

\item  $\mathcal{I}_{1}\cap \left( 2,+\infty \right) \neq \varnothing $ if
and only if $a_{0}\in \mathbb{R}$ and $b_{0}>\underline{b}\left( a_{0}\right) $;

\item  $\mathcal{I}_{2}\cap \left( 2,+\infty \right) \neq \varnothing $ for
every $a,b\in \mathbb{R}$;

\item  $\mathcal{I}_{1}\cap \left( 1,2\right) \neq \varnothing $ if and only
if $a_{0}\in \mathbb{R}$ and $b_{0}>\min \left\{ a_{0},-\frac{N-a_{0}}{2},-%
\frac{N+2}{2}\right\} $;

\item  $\mathcal{I}_{2}\cap \left( 1,2\right) \neq \varnothing $ if and only
if $a\in \mathbb{R}$ and $b<\max \left\{ a,-2\right\} $.
\end{itemize}

Our main existence result for super-linear nonlinearities is the following
theorem. A related result, concerning ground state solutions, will be given
in Section \ref{SEC: GS}.

\begin{theorem}
\label{THM: super general}Assume $\left( \mathbf{V}\right) ,\left( \mathbf{K}%
\right) $ with $a_{0},a,b\in \mathbb{R}$ and $b_{0}>\underline{b}\left(
a_{0}\right) $. Then Eq. (\ref{Eq.}) has a nonnegative nontrivial radial
solution for every continuous $f:\mathbb{R}\rightarrow \mathbb{R}$ satisfying:

\begin{itemize}
\item[$\left( \mathbf{f}_{3}\right) $]  
$\displaystyle \sup_{t>0}\,\frac{\left| f\left( t\right) \right| }{\min \left\{
t^{q_{1}-1},t^{q_{2}-1}\right\} }<+\infty $ for some $q_{1},q_{2}>2$ such
that $q_{1}\in \mathcal{I}_{1}$, $q_{2}\in \mathcal{I}_{2};$

\item[$\left( \mathbf{f}_{4}\right) $]  
$\exists \theta >2$ such that $0\leq \theta F\left( t\right) \leq f\left( t\right) t$ for all $t>0;$

\item[$\left( \mathbf{f}_{5}\right) $]  
$\exists t_{0}>0$ such that $F\left(t_{0}\right) >0.$
\end{itemize}

\noindent If $K\left( \left| \cdot \right| \right) \in L^{1}(\mathbb{R}^{N})$,
the same result holds with $\left( \mathbf{f}_{4}\right) $-$\left( \mathbf{f}%
_{5}\right) $ replaced by the weaker condition:

\begin{itemize}
\item[$\left( \mathbf{f}_{6}\right) $]  $\exists \theta >2$ and $\exists
t_{0}>0$ such that $0<\theta F\left( t\right) \leq f\left( t\right) t$ for
all $t\geq t_{0}.$
\end{itemize}
\end{theorem}

Observe that, as soon as we can take $q_{1}=q_{2}$ in $\left( \mathbf{f}%
_{3}\right) $, i.e., when $\mathcal{I}_{1}\cap \mathcal{I}_{2}\neq
\varnothing $, Theorem \ref{THM: super general} only requires that the
nonlinearity satisfies a single power growth condition, i.e., a condition of
the form 
\[
\left| f\left( t\right) \right| \leq \left( const.\right) t^{q-1}\qquad 
\textrm{with\quad }q\in \mathcal{I}_{1}\cap \mathcal{I}_{2}\cap \left(
2,+\infty \right) . 
\]
Indeed, such a condition is equivalent to $\left( \mathbf{f}_{3}\right) $,
because: (i) it obviously implies $\left( \mathbf{f}_{3}\right) $ with $%
q_{1}=q_{2}=q$; (ii) assuming for instance that $\left( \mathbf{f}%
_{3}\right) $ holds with $q_{1}\leq q_{2}$, one has $\min \left\{
t^{q_{1}-1},t^{q_{2}-1}\right\} \leq t^{q-1}$ for every $t>0$ and $q\in
\left[ q_{1},q_{2}\right] $, and one can find such a $q$ which also belongs
to $\mathcal{I}_{1}\cap \mathcal{I}_{2}\cap \left( 2,+\infty \right) $.

Moreover, assuming $b_{0}>\underline{b}\left( a_{0}\right) $, it is not
difficult to check that $\mathcal{I}_{1}\cap \mathcal{I}_{2}\neq \varnothing 
$ if and only if $\underline{q}<\overline{q}$ (where $\underline{q}$ and $%
\overline{q}$ are the exponents of Theorem \ref{THM: SWW}) and, in this
case, one has 
\[
\mathcal{I}_{1}\cap \mathcal{I}_{2}\cap \left( 2,+\infty \right) =(%
\underline{q},\overline{q}). 
\]
Therefore, the case $\mathcal{I}_{1}\cap \mathcal{I}_{2}\neq \varnothing $
is exactly the case in which the potentials behave compatibly at zero and
infinity, and, in such a case, since $\left( \mathbf{f}_{3}\right) $ becomes
equivalent to $\left( \mathbf{f}_{1}\right) $, the first part of Theorem \ref
{THM: super general} gives exactly Theorem \ref{THM: SWW} (up to the fact
that the pair $\left( \mathbf{f}_{4}\right) $-$\left( \mathbf{f}_{5}\right) $
is a slightly weaker condition than $\left( \mathbf{f}_{2}\right) $, which
is a rather technical generalization and it is not what we want to emphasize
here). Note that $\mathcal{I}_{1}$ and $\mathcal{I}_{2}$ only depend on $%
a_{0},b_{0}$ and $a,b$ respectively, so that $\mathcal{I}_{1}\cap \mathcal{I}%
_{2}\neq \varnothing $ means a link between the potential rates at zero and
infinity (precisely: $q_{**}\left( a,b\right) <q^{*}\left(
a_{0},b_{0}\right) $).

The case $\mathcal{I}_{1}\cap \mathcal{I}_{2}=\varnothing $, instead, is the
case without compatibility (and amounts to $\overline{q}\leq \underline{q}$%
), so that Theorem \ref{THM: SWW} does not apply and Theorem \ref{THM: super
general} is a new result (both in its first and second part), concerning
nonlinearities that satisfy the double power growth condition $\left( 
\mathbf{f}_{3}\right) $ with $q_{1}\neq q_{2}$. For the reader convenience,
we explain this case in the following corollary.

\begin{corollary}
\label{COR: double}Assume $\left( \mathbf{V}\right) ,\left( \mathbf{K}%
\right) $ with $a_{0}>-\left( 2N-2\right) $, $b_{0}>\min \left\{
a_{0},-2\right\} $ and one of the following alternatives: 
\begin{equation}
a\leq -2,\qquad b\geq \max \left\{ 2\frac{\left( N-2\right) b_{0}-\left(
N-1\right) \left( a_{0}+2\right) }{2N-2+a_{0}},b_{0}\right\} 
\label{COR:hp1}
\end{equation}
or 
\begin{equation}
b>a>-2,\qquad \frac{b-a}{2N-2+a}\geq \max \left\{ \frac{b_{0}-a_{0}}{%
2N-2+a_{0}},\frac{b_{0}+2}{2\left( N-2\right) }\right\} .  \label{COR:hp2}
\end{equation}
Then Eq. (\ref{Eq.}) has a nonnegative nontrivial radial solution for every
continuous $f:\mathbb{R}\rightarrow \mathbb{R}$ satisfying $\left( \mathbf{f}%
_{4}\right) $, $\left( \mathbf{f}_{5}\right) $ and 
\[
\sup_{t>0}\,\frac{\left| f\left( t\right) \right| }{\min \left\{
t^{q_{1}-1},t^{q_{2}-1}\right\} }<+\infty 
\]
for some 
\begin{equation}
2<q_{1}<\max \left\{ 2\frac{N+b_{0}}{N-2},2\frac{2N-2+2b_{0}-a_{0}}{%
2N-2+a_{0}}\right\} \quad \textrm{and}\quad q_{2}>\min \left\{ 2\frac{N+b}{N-2}%
,2\frac{2N-2+2b-a}{2N-2+a}\right\} .  \label{COR:th}
\end{equation}
If $K\left( \left| \cdot \right| \right) \in L^{1}(\mathbb{R}^{N})$, the same
result holds with $\left( \mathbf{f}_{4}\right) $-$\left( \mathbf{f}%
_{5}\right) $ replaced by $\left( \mathbf{f}_{6}\right) $.
\end{corollary}

The interested reader can check that Corollary \ref{COR: double} is exactly
the case of Theorem \ref{THM: super general} with $\mathcal{I}_{1}\cap 
\mathcal{I}_{2}=\varnothing $. We just observe that, under the assumptions
of the corollary, one explicitly has 
\begin{eqnarray}
\max \left\{ \frac{b_{0}-a_{0}}{2N-2+a_{0}},\frac{b_{0}+2}{2\left(
N-2\right) }\right\} &=&\left\{ 
\begin{array}{lll}
\frac{b_{0}-a_{0}}{2N-2+a_{0}} & ~\smallskip & \textrm{if }a_{0}<-2 \\ 
\frac{b_{0}+2}{2\left( N-2\right) } &  & \textrm{if }a_{0}\geq -2,
\end{array}
\right.  \nonumber \\
\max \left\{ 2\frac{\left( N-2\right) b_{0}-\left( N-1\right) \left(
a_{0}+2\right) }{2N-2+a_{0}},b_{0}\right\} &=&\left\{ 
\begin{array}{lll}
2\frac{\left( N-2\right) b_{0}-\left( N-1\right) \left( a_{0}+2\right) }{%
2N-2+a_{0}} & ~\smallskip & \textrm{if }a_{0}<-2 \\ 
b_{0} &  & \textrm{if }a_{0}\geq -2,
\end{array}
\right.  \nonumber \\
\max \left\{ 2\frac{N+b_{0}}{N-2},2\frac{2N-2+2b_{0}-a_{0}}{2N-2+a_{0}}%
\right\} &=&\left\{ 
\begin{array}{lll}
2\frac{2N-2+2b_{0}-a_{0}}{2N-2+a_{0}} & ~\smallskip & \textrm{if }a_{0}<-2 \\ 
2\frac{N+b_{0}}{N-2} &  & \textrm{if }a_{0}\geq -2,
\end{array}
\right.  \label{Cup} \\
\min \left\{ 2\frac{N+b}{N-2},2\frac{2N-2+2b-a}{2N-2+a}\right\} &=&\left\{ 
\begin{array}{lll}
2\frac{N+b}{N-2} & ~\smallskip & \textrm{if (\ref{COR:hp1}) holds} \\ 
2\frac{2N-2+2b-a}{2N-2+a} &  & \textrm{if (\ref{COR:hp2}) holds,}
\end{array}
\right.  \label{Cdwn}
\end{eqnarray}
where (\ref{Cup}) and (\ref{Cdwn}), which are the right hand sides of
inequalities (\ref{COR:th}), respectively coincide with $q^{*}\left(
a_{0},b_{0}\right) $\ and $q_{**}\left( a,b\right) $, or equivalently with $%
\overline{q}\left( a_{0},b_{0}\right) $ and $\underline{q}\left(
a,b,a_{0},b_{0}\right) $, and satisfy $q^{*}\left( a_{0},b_{0}\right) \leq
q_{**}\left( a,b\right) $.\smallskip

As far as sub-linear nonlinearities are concerned, we will prove the
following result.

\begin{theorem}
\label{THM: sub general}Assume $\left( \mathbf{V}\right) ,\left( \mathbf{K}%
\right) $ with $a_{0},a\in \mathbb{R}$, $b_{0}>\min \left\{ a_{0},-\frac{N-a_{0}%
}{2},-\frac{N+2}{2}\right\} $ and $b<\max \left\{ a,-2\right\} $. Then Eq. (%
\ref{Eq.}) has a nonnegative nontrivial radial solution for every continuous 
$f:\mathbb{R}\rightarrow \mathbb{R}$ satisfying:

\begin{itemize}
\item[$\left( \mathbf{f}_{7}\right) $]  
$\displaystyle \sup_{t>0}\,\frac{\left| f\left( t\right) \right| }{\min \left\{
t^{q_{1}-1},t^{q_{2}-1}\right\} }<+\infty $ for some $q_{1},q_{2}<2$ such
that $q_{1}\in \mathcal{I}_{1}$, $q_{2}\in \mathcal{I}_{2};$

\item[$\left( \mathbf{f}_{8}\right) $]  
$\exists \theta <2$ such that $\displaystyle\liminf_{t\rightarrow 0^{+}}\frac{F\left( t\right) }{t^{\theta }}>0.$
\end{itemize}
\end{theorem}

In contrast with the case of Theorem \ref{THM: super general} with respect
to Theorem \ref{THM: SWW}, Theorem \ref{THM: sub general} extends Theorem 
\ref{THM: ST} in many directions (other than the obvious fact that Theorem 
\ref{THM: sub general} concerns more general nonlinearities than the pure
power $f\left( u\right) =\left| u\right| ^{q-2}u$). Such improvements are
described by the following remarks, where the set 
\[
\mathcal{P}=\left\{ \left( a,b,a_{0},b_{0}\right) :b_{0}>\min \left\{ a_{0},-%
\frac{N-a_{0}}{2},-\frac{N+2}{2}\right\} ,\,b<\max \left\{ a,-2\right\} ,\,%
\mathcal{I}_{1}\cap \mathcal{I}_{2}\neq \varnothing \right\} 
\]
is used. Notice that $\mathcal{P}$ is the set of the potential rates $%
a,b,a_{0},b_{0}$ for which Theorem \ref{THM: sub general} concerns single
power nonlinearities, in the sense that, if $\left( a,b,a_{0},b_{0}\right)
\in \mathcal{P}$, then $\left( \mathbf{f}_{7}\right) $ is equivalent to 
\[
\left| f\left( t\right) \right| \leq \left( const.\right) t^{q-1}\textrm{%
\qquad with\quad }q\in \mathcal{I}_{1}\cap \mathcal{I}_{2}\cap \left(
1,2\right) 
\]
(cf. the discussion after Theorem \ref{THM: super general}).

\begin{itemize}
\item  \emph{The set }$\mathcal{P}$\emph{\ is strictly larger than the set
of the potential rates for which Theorem \ref{THM: ST} applies, i.e., the set%
}\textsl{\ } 
\[
\mathcal{P}_{1}=\left\{ \left( a,b,a_{0},b_{0}\right) :\left( a,b\right) \in 
\mathcal{A}_{1}\cup ...\cup \mathcal{A}_{5},\,\left( a_{0},b_{0}\right) \in 
\mathcal{B}_{1}\cup ...\cup \mathcal{B}_{6},\,\underline{\underline{q}}<%
\overline{\overline{q}}\,\right\} .
\]
For example, if $\left( a,b\right) \in \mathcal{A}_{2}$ and $\left(
a_{0},b_{0}\right) \in \mathcal{B}_{1}$ with $a_{0}<-\left( 2N-2\right) $
and $b\geq b_{0}$, then one has 
\begin{equation}
\underline{\underline{q}}=2\frac{N+b}{N-2},\quad \overline{\overline{q}}=2%
\frac{N+b_{0}}{N-2},\quad q_{*}=1,\quad q^{*}=+\infty ,\quad q_{**}=1
\label{EX 1}
\end{equation}
and therefore $\mathcal{I}_{1}\cap \mathcal{I}_{2}\cap \left( 1,2\right)
=\left( 1,2\right) $ and $\underline{\underline{q}}\geq \overline{\overline{q%
}}$, so that $\left( a,b,a_{0},b_{0}\right) \in \mathcal{P}$ but $\left(
a,b,a_{0},b_{0}\right) \notin \mathcal{P}_{1}$.

Other examples are given by those points $\left( a,b,a_{0},b_{0}\right) \in 
\mathcal{P}$ for which at least one of the exponents $\underline{\underline{q%
}}$ and $\overline{\overline{q}}$ is undefined, i.e., $\left(
a,b,a_{0},b_{0}\right) \notin \left( \mathcal{A}_{1}\cup ...\cup \mathcal{A}%
_{5}\right) \times \left( \mathcal{B}_{1}\cup ...\cup \mathcal{B}_{6}\right) 
$. For instance, if 
\[
\frac{a-2}{2}\leq b<a,\quad a_{0}<-\left( 2N-2\right) \quad \textrm{and\quad }%
b_{0}>-2,
\]
then both $\underline{\underline{q}}$ and $\overline{\overline{q}}$ are
undefined, while we get $q_{*}=1$, $q^{*}=+\infty $, $q_{**}=2\frac{N+b}{N+a}
$ and therefore $\mathcal{I}_{1}\cap \mathcal{I}_{2}\cap \left( 1,2\right)
=(2\frac{N+b}{N+a},2)$.

This means that Theorem \ref{THM: sub general} gives existence results to
Eq. (\ref{Eq.}) with power type nonlinearities (even with $f\left( u\right)
=\left| u\right| ^{q-2}u$) for more potentials than the ones allowed by
Theorem \ref{THM: ST}, and exactly for those potentials $V,K$ satisfying $%
\left( \mathbf{V}\right) ,\left( \mathbf{K}\right) $ with $\left(
a,b,a_{0},b_{0}\right) \in \mathcal{P}\setminus \mathcal{P}_{1}$. The
explicit description of the whole set $\mathcal{P}\setminus \mathcal{P}_{1}$
is left to the interested reader. We just observe that the above examples
show that $\mathcal{P}\setminus \mathcal{P}_{1}$ contains both points for
which $\underline{\underline{q}}\geq \overline{\overline{q}}$ and points for
which $\underline{\underline{q}}$ and $\overline{\overline{q}}$ are
undefined.

\item  \emph{If }$\left( a,b,a_{0},b_{0}\right) \in \mathcal{P}_{1}$\emph{,
the interval }$\mathcal{I}_{1}\cap \mathcal{I}_{2}\cap \left( 1,2\right) $%
\emph{\ can be strictly larger than }$(\underline{\underline{q}},\overline{%
\overline{q}})$\emph{\ (depending on }$a,b,a_{0},b_{0}$\emph{)}. For
example, if we take $\left( a,b\right) \in \mathcal{A}_{1}$ and $\left(
a_{0},b_{0}\right) \in \mathcal{B}_{2}$ with $b<b_{0}$, then 
\[
\underline{\underline{q}}=2\frac{N+b}{N-2},\quad \overline{\overline{q}}=2%
\frac{N+b_{0}}{N-2},\quad q_{*}=1,\quad q^{*}=2\frac{N+b_{0}}{N-2},\quad
q_{**}=2\frac{N+b}{N-2}
\]
and therefore 
\[
\mathcal{I}_{1}\cap \mathcal{I}_{2}\cap \left( 1,2\right) =\left( 1,2\frac{%
N+b_{0}}{N-2}\right) \cap \left( 2\frac{N+b}{N-2},+\infty \right) \cap
\left( 1,2\right) =\left( 2\frac{N+b}{N-2},2\frac{N+b_{0}}{N-2}\right) =(%
\underline{\underline{q}},\overline{\overline{q}}).
\]
But if we take $\left( a,b\right) \in \mathcal{A}_{2}$ and $\left(
a_{0},b_{0}\right) \in \mathcal{B}_{1}$ with $a_{0}<-\left( 2N-2\right) $
and $-\frac{N-2}{2}<b<b_{0}<-2$, then we have (\ref{EX 1}) as before and
therefore $\mathcal{I}_{1}\cap \mathcal{I}_{2}\cap \left( 1,2\right) =\left(
1,2\right) $ strictly contains $(\underline{\underline{q}},\overline{%
\overline{q}})=(2\frac{N+b}{N-2},2\frac{N+b_{0}}{N-2})$.

This means that there are potentials for which Theorem \ref{THM: ST} applies
but Theorem \ref{THM: sub general} gives a wider range of power type
nonlinearities for which Eq. (\ref{Eq.}) admits solutions, and exactly those
potentials $V,K$ satisfying $\left( \mathbf{V}\right) ,\left( \mathbf{K}%
\right) $ with $a,b,a_{0},b_{0}$ such that $\mathcal{I}_{1}\cap \mathcal{I}%
_{2}\cap \left( 1,2\right) \setminus (\underline{\underline{q}},\overline{%
\overline{q}})$ is nonempty. We leave to the interested reader the explicit
description of the set $\mathcal{I}_{1}\cap \mathcal{I}_{2}\cap \left(
1,2\right) \setminus (\underline{\underline{q}},\overline{\overline{q}})$,
as $a,b,a_{0},b_{0}$ vary.

\item  \emph{If }$a,b,a_{0},b_{0}$\emph{\ are such that Theorem \ref{THM:
sub general} applies with }$\mathcal{I}_{1}\cap \mathcal{I}_{2}=\varnothing $%
\emph{, then one can check that Theorem \ref{THM: ST} does not apply}, so
that Theorem \ref{THM: sub general} gives new existence results, concerning
double power nonlinearities that satisfy the growth condition $\left( 
\mathbf{f}_{7}\right) $ with $q_{1}\neq q_{2}$. The explicit description of
the set of the potential rates $a_{0},a\in \mathbb{R}$, $b_{0}>\min \left\{
a_{0},-\frac{N-a_{0}}{2},-\frac{N+2}{2}\right\} $ and $b<\max \left\{
a,-2\right\} $ for which $\mathcal{I}_{1}\cap \mathcal{I}_{2}=\varnothing $
is left to the interested reader.
\end{itemize}

Theorems \ref{THM: super general} and \ref{THM: sub general} will be proved
in Section \ref{SEC: proofs} by variational methods, as a first application
of the compactness results of \cite{BGR-p1}. Other applications will be
given in \cite{BGR-p2}, where Eq. (\ref{Eq.}) will be studied with more
general potentials (not necessarily continuous and possibly not satisfying
power type estimates at the origin and at infinity) and nonlinearities
(e.g., the presence of an additional forcing term is considered), also
dealing with the case of bounded and exterior domains. A version of Theorem 
\ref{THM: super general} without the \textit{Ambrosetti-Rabinowitz condition}
$\left( \mathbf{f}_{4}\right) $ will be given in \cite{GR-boundedPS}.

The proof of Theorems \ref{THM: super general} and \ref{THM: sub general}
will be achieved here by the same techniques used in \cite{Su-Wang-Will
p,SuTian12} for proving Theorems \ref{THM: SWW} and \ref{THM: ST}, namely,
respectively, the application of the Mountain Pass Theorem and the global
minimization on $H_{V,\mathrm{r}}^{1}$ of the Euler functional associated to
the equation. The main difference between our arguments and the ones of \cite
{Su-Wang-Will p,SuTian12} is that the single power growth assumption
required on the nonlinearity in Theorems \ref{THM: SWW} and \ref{THM: ST}
only allows to exploit the compact embedding of $H_{V,\mathrm{r}}^{1}$ into
the weighted Lebesgue space $L_{K}^{q}:=L^{q}(\mathbb{R}^{N},K\left( \left|
x\right| \right) dx)$, while the double power growth assumptions $\left( 
\mathbf{f}_{3}\right) $ and $\left( \mathbf{f}_{7}\right) $ allow us to use
the more general compact embedding \cite{BGR-p1} of $H_{V,\mathrm{r}}^{1}$
into the sum space $L_{K}^{q_{1}}+L_{K}^{q_{2}}$ (see Section \ref{SEC:
proofs} for some recallings on such a space). The fact that, for $%
q_{1}=q_{2}=q$, the space $L_{K}^{q_{1}}+L_{K}^{q_{2}}$ becomes $L_{K}^{q}$
and the double power growth assumption becomes the single power one reflects
on a technical level the already discussed fact that Theorems \ref{THM:
super general} and \ref{THM: sub general} contain Theorems \ref{THM: SWW}
and \ref{THM: ST} and extend them to a wider class of potentials, avoiding
any compatibility requirement between their behaviours at the origin and at
infinity.

For all the considerations expounded in this introduction, we believe that
the double power growth assumption and the related sum space $%
L_{K}^{q_{1}}+L_{K}^{q_{2}}$ are the ``right'' tools for studying problems
like (\ref{Eq.}), i.e., problems on the whole space in which some weights
are present and both their behaviours at zero and at infinity affect the
solutions (for a different use of the sum of Lebesgue spaces in nonlinear
problems, see \cite{Azz-D'Av-Pomp14,Zhang-Fu15}).\smallskip

We conclude the section with some remarks and examples of nonlinearities
satisfying our assumptions.

\begin{remark}
\label{RMK: finale}\quad 

\begin{enumerate}
\item  Under the same assumptions of Theorem \ref{THM: super general}, if $f$
is also odd and satisfies 
\begin{equation}
\inf_{t>0}\,\frac{f\left( t\right) }{\min \left\{
t^{q_{1}-1},t^{q_{2}-1}\right\} }>0  \label{mult}
\end{equation}
(with the same exponents of hypothesis $\left( \mathbf{f}_{3}\right) $),
then Eq. (\ref{Eq.}) has infinitely many radial solutions. Similarly, if the
same assumptions of Theorem \ref{THM: sub general} hold and if $f$ is also
odd, then Eq. (\ref{Eq.}) has infinitely many radial solutions. These
results rely on the variational theory of symmetric functionals and we refer
the reader to the analogous results of \cite{BGR-p2} for a detailed proof.

\item  \label{RMK:symm-crit}The solutions found in both Theorems \ref{THM:
super general} and \ref{THM: sub general} also satisfy (\ref{weak solution})
for all $h\in H_{V}^{1}$, since, under the hypotheses of the theorems, the
symmetric criticality type results of \cite{BGR-p2} apply.

\item  The continuity of $V$ and $K$ is not essential to Theorems \ref{THM:
super general} and \ref{THM: sub general}, and may be replaced by some
weaker integrability assumptions. We refer the interested reader again to 
\cite{BGR-p2} for a generalization of Theorems \ref{THM: super general} and 
\ref{THM: sub general} in this direction.
\end{enumerate}
\end{remark}

\begin{example}
\label{EX:f}The more obvious function with a double power growth is 
$f\left(t\right) =\min \left\{ \left| t\right| ^{q_{1}-2}t,\left| t\right|^{q_{2}-2}t\right\} $, 
which also satisfies $\left( \mathbf{f}_{4}\right) $
(with $\theta =\min \left\{ q_{1},q_{2}\right\} $) if $q_{1},q_{2}>2$, and 
$\left( \mathbf{f}_{8}\right) $ (with $\theta =\max \left\{q_{1},q_{2}\right\} $) if $q_{1},q_{2}<2$.
Another model example is 
\[
f\left( t\right) =\frac{\left| t\right| ^{q_{2}-2}t}{1+\left| t\right|
^{q_{2}-q_{1}}}\quad \textrm{with }q_{1}\leq q_{2},
\]
for which $\left( \mathbf{f}_{4}\right) $ holds (with $\theta =q_{1}$) if $%
q_{1}>2$ and $\left( \mathbf{f}_{8}\right) $ holds (with $\theta =q_{2}$) if 
$q_{2}<2$.
Note that both these functions are odd and also satisfy (\ref{mult}). 
Moreover, both of them become $f\left( t\right) =\left| t\right| ^{q-2}t$ if $q_{1}=q_{2}=q$.
Other examples of nonlinearities satisfying $%
\sup_{t>0}\,\left| f\left( t\right) \right| /\min \left\{
t^{q_{1}-1},t^{q_{2}-1}\right\} <+\infty $ are
\[
f\left( t\right) =\frac{\left| t\right| ^{q_{1}+q-1}-\left| t\right|
^{q_{2}-1}}{1+\left| t\right| ^{q}},\quad f\left( t\right) =\frac{\left|
t\right| ^{q_{2}-1+\varepsilon }}{1+\left| t\right|
^{q_{2}-q_{1}+2\varepsilon }}\ln \left| t\right| 
\]
(the latter extended at $0$ by continuity) with $1<q_{1}\leq q_{2}<q_{1}+q$
and $\varepsilon >0$, which do not satisfy $\left( \mathbf{f}_{4}\right) $
or $\left( \mathbf{f}_{8}\right) $, but satisfy $\left( \mathbf{f}%
_{6}\right) $ if $q_{1}>2$ and $\varepsilon $ is small enough (precisely: $%
\varepsilon <q_{1}-2$).
\end{example}

\section{Proof of Theorems \ref{THM: super general} and \ref{THM: sub
general} \label{SEC: proofs}}

{\allowdisplaybreaks

Let $N\geq 3$ and let $V,K$ be as in $\left( \mathbf{V}\right) ,\left( 
\mathbf{K}\right) $ with $a_{0},a,b\in \mathbb{R}$ and $b_{0}>b_{*}\left(
a_{0}\right) $. Recall the definition (\ref{H^1_V,r}) of $H_{V,\mathrm{r}%
}^{1}$, which is a Hilbert space with respect the following inner product
and related norm: 
\begin{equation}
\left( u\mid v\right) :=\int_{\mathbb{R}^{N}}\nabla u\cdot \nabla v\,dx+\int_{%
\mathbb{R}^{N}}V\left( \left| x\right| \right) uv\,dx,\quad \left\| u\right\|
:=\left( \int_{\mathbb{R}^{N}}\left| \nabla u\right| ^{2}dx+\int_{\mathbb{R}%
^{N}}V\left( \left| x\right| \right) u^{2}dx\right) ^{1/2}.  \label{h struct}
\end{equation}

Denote by $L_{K}^{q}(\mathbb{R}^{N}):=L^{q}(\mathbb{R}^{N},K\left( \left| x\right|
\right) dx)$ the usual Lebesgue space with respect to the measure $K\left(
\left| x\right| \right) dx$ ($dx$ stands for the Lebesgue measure on $\mathbb{R}%
^{N}$) and consider the sum space 
\[
L_{K}^{q_{1}}+L_{K}^{q_{2}}:=\left\{ u_{1}+u_{2}:u_{1}\in
L_{K}^{q_{1}}\left( \mathbb{R}^{N}\right) ,\,u_{2}\in L_{K}^{q_{2}}\left( \mathbb{R%
}^{N}\right) \right\} ,\quad 1<q_{i}<\infty . 
\]
From \cite{BPR}, we recall that such a space is a Banach space with respect
to the norm 
\[
\left\| u\right\| _{L_{K}^{q_{1}}+L_{K}^{q_{2}}}:=\inf_{u_{1}+u_{2}=u}\max
\left\{ \left\| u_{1}\right\| _{L_{K}^{q_{1}}(\mathbb{R}^{N})},\left\|
u_{2}\right\| _{L_{K}^{q_{2}}(\mathbb{R}^{N})}\right\} 
\]
and can be characterized as the set of the measurable mappings $u:\mathbb{R}%
^{N}\rightarrow \mathbb{R}$ for which there exists a measurable set $E\subseteq 
\mathbb{R}^{N}$ such that $u\in L_{K}^{q_{1}}\left( E\right) \cap
L_{K}^{q_{2}}\left( E^{c}\right) $.

Recall the definitions (\ref{I_1:=}) and (\ref{I_2:=}) of the intervals $%
\mathcal{I}_{1}=\mathcal{I}_{1}\left( a_{0},b_{0}\right) $ and $\mathcal{I}%
_{2}=\mathcal{I}_{2}\left( a,b\right) $.

\begin{lemma}
\label{LEM: S->0}For every $q_{1}\in \mathcal{I}_{1}$ and $q_{2}\in \mathcal{%
I}_{2}$ one has $\displaystyle\lim_{R\rightarrow 0^{+}}\mathcal{S}_{1}\left( R\right) =\lim_{R\rightarrow
+\infty }\mathcal{S}_{2}\left( R\right) =0$, where 
\[
\mathcal{S}_{1}\left( R\right) :=\sup_{u\in H_{V,\mathrm{r}}^{1},\,\left\|
u\right\| =1}\int_{B_{R}}K\left( \left| x\right| \right) \left| u\right|
^{q_{1}}dx,\quad \mathcal{S}_{2}\left( R\right) :=\sup_{u\in H_{V,\mathrm{r}%
}^{1},\,\left\| u\right\| =1}\int_{\mathbb{R}^{N}\setminus B_{R}}K\left( \left|
x\right| \right) \left| u\right| ^{q_{2}}dx.
\]
\end{lemma}

\proof%
It follows from the results of \cite{BGR-p1}, and precisely from Theorem 4
(apply with $\alpha _{\infty }=b$, $\beta _{\infty }=0$, $\gamma _{\infty
}=-a$ if $a>-2$ and $\gamma _{\infty }=2$ if $a\leq -2$) and Theorem 5
(apply with $\alpha _{0}=b_{0}$, $\beta _{0}=0$, $\gamma _{0}=-a_{0}$ if $%
a_{0}<-2$ and $\gamma _{0}=2$ if $a_{0}\geq -2$).%
\endproof

\begin{lemma}
\label{LEM: comp embedd}The space $H_{V,\mathrm{r}}^{1}$ is compactly
embedded into $L_{K}^{q_{1}}+L_{K}^{q_{2}}$ for every $q_{1}\in \mathcal{I}%
_{1}$ and $q_{2}\in \mathcal{I}_{2}$.
\end{lemma}

\proof%
It readily follows from Lemma \ref{LEM: S->0} above and Theorem 1 of \cite
{BGR-p1}.%
\endproof
\bigskip

Now assume that $f:\mathbb{R}\rightarrow \mathbb{R}$ is a continuous function for
which there exist $q_{1}\in \mathcal{I}_{1}$, $q_{2}\in \mathcal{I}_{2}$ and 
$M>0$ such that 
\begin{equation}
\left| f\left( t\right) \right| \leq M\min \left\{ \left| t\right|
^{q_{1}-1},\left| t\right| ^{q_{2}-1}\right\} \qquad \textrm{for all }t\in 
\mathbb{R}.  \label{|f| <}
\end{equation}
Set $F\left( t\right) :=\int_{0}^{t}f\left( s\right) ds$ and define the
functional 
\begin{equation}
I\left( u\right) :=\frac{1}{2}\left\| u\right\| ^{2}-\int_{\mathbb{R}%
^{N}}K\left( \left| x\right| \right) F\left( u\right) dx\qquad \textrm{for
every }u\in H_{V,\mathrm{r}}^{1}.  \label{I:=}
\end{equation}

\begin{lemma}
\label{LEM: C^1}$I$ is a $C^{1}$ functional on $H_{V,\mathrm{r}}^{1}$ and
its Fr\'{e}chet derivative $I^{\prime }\left( u\right) $ at any $u\in H_{V,%
\mathrm{r}}^{1}$ is given by 
\[
I^{\prime }\left( u\right) h=\int_{\mathbb{R}^{N}}\nabla u\cdot \nabla
h\,dx+\int_{\mathbb{R}^{N}}V\left( \left| x\right| \right) uh\,dx-\int_{\mathbb{R}%
^{N}}K\left( \left| x\right| \right) f\left( u\right) h\,dx,\quad \forall
h\in H_{V,\mathrm{r}}^{1}.
\]
\end{lemma}

\proof%
It follows from Lemma \ref{LEM: comp embedd} above and the results of \cite
{BPR} about Nemytski\u{\i} operators on the sum of Lebesgue spaces. Indeed,
by \cite[Proposition 3.8]{BPR}, condition (\ref{|f| <}) implies that the
functional 
\[
u\in L_{K}^{q_{1}}+L_{K}^{q_{2}}\mapsto \int_{\mathbb{R}^{N}}K\left( \left|
x\right| \right) F\left( u\right) dx 
\]
is of class $C^{1}$ with Fr\'{e}chet derivative at any $u\in
L_{K}^{q_{1}}+L_{K}^{q_{2}}$ is given by 
\[
h\in L_{K}^{q_{1}}+L_{K}^{q_{2}}\mapsto \int_{\mathbb{R}^{N}}K\left( \left|
x\right| \right) f\left( u\right) h\,dx. 
\]
The result then ensues by the continuous embedding $H_{V,\mathrm{r}%
}^{1}\hookrightarrow L_{K}^{q_{1}}+L_{K}^{q_{2}}$ given by Lemma \ref{LEM:
comp embedd}.%
\endproof%
\bigskip

By Lemma \ref{LEM: C^1}, the problem of finding radial solutions to Eq. (\ref
{Eq.}) clearly reduces to the problem of finding critical points of $I:H_{V,%
\mathrm{r}}^{1}\rightarrow \mathbb{R}$.

For future reference, we observe here that, by condition (\ref{|f| <}),
there exists $\tilde{M}>0$ such that 
\begin{equation}
\left| F\left( t\right) \right| \leq \tilde{M}\min \left\{ \left| t\right|
^{q_{1}},\left| t\right| ^{q_{2}}\right\} \qquad \textrm{for all }t\in \mathbb{R}.
\label{|F| <}
\end{equation}

\begin{lemma}
\label{LEM: I(u)>}There exist two constants $c_{1},c_{2}>0$ such that 
\begin{equation}
I\left( u\right) \geq \frac{1}{2}\left\| u\right\| ^{2}-c_{1}\left\|
u\right\| ^{q_{1}}-c_{2}\left\| u\right\| ^{q_{2}}\qquad \textrm{for all }u\in
H_{V,\mathrm{r}}^{1}.  \label{I(u)>}
\end{equation}
\end{lemma}

\proof%
By Lemma \ref{LEM: S->0}, fix $R_{2}>R_{1}>0$ such that $\mathcal{S}%
_{1}\left( R_{1}\right) ,\mathcal{S}_{2}\left( R_{2}\right) <1$. Then, by 
\cite[Lemma 1]{BGR-p1} and the continuous embedding $H_{V,\mathrm{r}%
}^{1}\hookrightarrow D^{1,2}(\mathbb{R}^{N})\hookrightarrow L_{\mathrm{loc}%
}^{2}(\mathbb{R}^{N})$, there exists a constant $c_{R_{1},R_{2}}>0$ such that 
\[
\int_{B_{R_{2}}\setminus B_{R_{1}}}K\left( \left| x\right| \right) \left|
u\right| ^{q_{1}}dx\leq c_{R_{1},R_{2}}\left\| u\right\| ^{q_{1}}\qquad 
\textrm{for all }u\in H_{V,\mathrm{r}}^{1}. 
\]
Therefore, by (\ref{|F| <}) and the definitions of $\mathcal{S}_{1}$ and $%
\mathcal{S}_{2}$, for every $u\in H_{V,\mathrm{r}}^{1}$ we get 
\begin{eqnarray*}
\left| \int_{\mathbb{R}^{N}}K\left( \left| x\right| \right) F\left( u\right)
dx\right| &\leq &\tilde{M}\int_{\mathbb{R}^{N}}K\left( \left| x\right| \right)
\min \left\{ \left| u\right| ^{q_{1}},\left| u\right| ^{q_{2}}\right\} dx \\
&\leq &\tilde{M}\left( \int_{B_{R_{1}}}K\left( \left| x\right| \right)
\left| u\right| ^{q_{1}}dx+\int_{B_{R_{2}}^{c}}K\left( \left| x\right|
\right) \left| u\right| ^{q_{2}}dx+\int_{B_{R_{2}}\setminus
B_{R_{1}}}K\left( \left| x\right| \right) \left| u\right| ^{q_{1}}dx\right)
\\
&\leq &\tilde{M}\left( \left\| u\right\| ^{q_{1}}\int_{B_{R_{1}}}K\left(
\left| x\right| \right) \frac{\left| u\right| ^{q_{1}}}{\left\| u\right\|
^{q_{1}}}dx+\left\| u\right\| ^{q_{2}}\int_{B_{R_{2}}^{c}}K\left( \left|
x\right| \right) \frac{\left| u\right| ^{q_{2}}}{\left\| u\right\| ^{q_{2}}}%
dx+c_{R_{1},R_{2}}\left\| u\right\| ^{q_{1}}\right) \\
&\leq &\tilde{M}\left( \left\| u\right\| ^{q_{1}}\mathcal{S}_{1}\left(
R_{1}\right) +\left\| u\right\| ^{q_{2}}\mathcal{S}_{2}\left( R_{2}\right)
+c_{R_{1},R_{2}}\left\| u\right\| ^{q_{1}}\right) .
\end{eqnarray*}
This yields (\ref{I(u)>}).%
\endproof%
\bigskip

\begin{lemma}
\label{LEM:PS}Assume $f\left( t\right) =0$ for all $t<0$. If $f$ satisfies $%
\left( \mathbf{f}_{4}\right) $, or $K\left( \left| \cdot \right| \right) \in
L^{1}(\mathbb{R}^{N})$ and $f$ satisfies $\left( \mathbf{f}_{6}\right) $, then
the functional $I:H_{V,\mathrm{r}}^{1}\rightarrow \mathbb{R}$ satisfies the
Palais-Smale condition.
\end{lemma}

\proof%
Let $\left\{ u_{n}\right\} $ be a sequence in $H_{V,\mathrm{r}}^{1}$ such
that $\left\{ I\left( u_{n}\right) \right\} $ is bounded and $I^{\prime
}\left( u_{n}\right) \rightarrow 0$ in the dual space of $H_{V,\mathrm{r}%
}^{1}$. Hence 
\[
\frac{1}{2}\left\| u_{n}\right\| ^{2}-\int_{\mathbb{R}^{N}}K\left( \left|
x\right| \right) F\left( u_{n}\right) dx=O\left( 1\right) \quad \textrm{and}%
\quad \left\| u_{n}\right\| ^{2}-\int_{\mathbb{R}^{N}}K\left( \left| x\right|
\right) f\left( u_{n}\right) u_{n}dx=o\left( 1\right) \left\| u_{n}\right\|
. 
\]
If $f$ satisfies $\left( \mathbf{f}_{4}\right) $, then we have $\theta
F\left( t\right) \leq f\left( t\right) t$ for all $t\in \mathbb{R}$ (because $%
f\left( t\right) =0$ for $t<0$) and therefore we get 
\[
\frac{1}{2}\left\| u_{n}\right\| ^{2}+O\left( 1\right) =\int_{\mathbb{R}%
^{N}}K\left( \left| x\right| \right) F\left( u_{n}\right) dx\leq \frac{1}{%
\theta }\int_{\mathbb{R}^{N}}K\left( \left| x\right| \right) f\left(
u_{n}\right) u_{n}dx=\frac{1}{\theta }\left\| u_{n}\right\| ^{2}+o\left(
1\right) \left\| u_{n}\right\| , 
\]
which implies that $\left\{ \left\| u_{n}\right\| \right\} $ is bounded,
since $\theta >2$. If $K\left( \left| \cdot \right| \right) \in L^{1}(\mathbb{R}%
^{N})$ and $f$ satisfies $\left( \mathbf{f}_{6}\right) $, then we have $%
\theta F\left( t\right) \leq f\left( t\right) t$ for all $\left| t\right|
\geq t_{0}$ (because $f\left( t\right) =0$ for $t<0$) and 
\begin{eqnarray*}
\int_{\left\{ \left| u_{n}\right| \geq t_{0}\right\} }K\left( \left|
x\right| \right) f\left( u_{n}\right) u_{n}dx &=&\int_{\mathbb{R}^{N}}K\left(
\left| x\right| \right) f\left( u_{n}\right) u_{n}dx-\int_{\left\{ \left|
u_{n}\right| <t_{0}\right\} }K\left( \left| x\right| \right) f\left(
u_{n}\right) u_{n}dx \\
&\leq &\int_{\mathbb{R}^{N}}K\left( \left| x\right| \right) f\left(
u_{n}\right) u_{n}dx+\int_{\left\{ \left| u_{n}\right| <t_{0}\right\}
}K\left( \left| x\right| \right) \left| f\left( u_{n}\right) u_{n}\right| dx
\\
&\leq &\int_{\mathbb{R}^{N}}K\left( \left| x\right| \right) f\left(
u_{n}\right) u_{n}dx+M\int_{\left\{ \left| u_{n}\right| <t_{0}\right\}
}K\left( \left| x\right| \right) \min \left\{ \left| u_{n}\right|
^{q_{1}},\left| u_{n}\right| ^{q_{2}}\right\} dx \\
&\leq &\int_{\mathbb{R}^{N}}K\left( \left| x\right| \right) f\left(
u_{n}\right) u_{n}dx+M\min \left\{ t_{0}^{q_{1}},t_{0}^{q_{2}}\right\}
\int_{\left\{ \left| u_{n}\right| <t_{0}\right\} }K\left( \left| x\right|
\right) dx \\
&\leq &\int_{\mathbb{R}^{N}}K\left( \left| x\right| \right) f\left(
u_{n}\right) u_{n}dx+M\min \left\{ t_{0}^{q_{1}},t_{0}^{q_{2}}\right\}
\left\| K\right\| _{L^{1}(\mathbb{R}^{N})},
\end{eqnarray*}
so that, by (\ref{|F| <}), we get 
\begin{eqnarray*}
\frac{1}{2}\left\| u_{n}\right\| ^{2}+O\left( 1\right) &=&\int_{\mathbb{R}%
^{N}}K\left( \left| x\right| \right) F\left( u_{n}\right) dx=\int_{\left\{
\left| u_{n}\right| <t_{0}\right\} }K\left( \left| x\right| \right) F\left(
u_{n}\right) dx+\int_{\left\{ \left| u_{n}\right| \geq t_{0}\right\}
}K\left( \left| x\right| \right) F\left( u_{n}\right) dx \\
&\leq &\tilde{M}\int_{\left\{ \left| u_{n}\right| <t_{0}\right\} }K\left(
\left| x\right| \right) \min \left\{ \left| u_{n}\right| ^{q_{1}},\left|
u_{n}\right| ^{q_{2}}\right\} dx+\frac{1}{\theta }\int_{\left\{ \left|
u_{n}\right| \geq t_{0}\right\} }K\left( \left| x\right| \right) f\left(
u_{n}\right) u_{n}dx \\
&\leq &\tilde{M}\min \left\{ t_{0}^{q_{1}},t_{0}^{q_{2}}\right\} \left\|
K\right\| _{L^{1}(\mathbb{R}^{N})}+\frac{1}{\theta }\int_{\mathbb{R}^{N}}K\left(
\left| x\right| \right) f\left( u_{n}\right) u_{n}dx+\frac{M}{\theta }\min
\left\{ t_{0}^{q_{1}},t_{0}^{q_{2}}\right\} \left\| K\right\| _{L^{1}(\mathbb{R}%
^{N})} \\
&=&\left( \tilde{M}+\frac{M}{\theta }\right) \min \left\{
t_{0}^{q_{1}},t_{0}^{q_{2}}\right\} \left\| K\right\| _{L^{1}(\mathbb{R}^{N})}+%
\frac{1}{\theta }\left\| u_{n}\right\| ^{2}+o\left( 1\right) \left\|
u_{n}\right\| .
\end{eqnarray*}
This yields again that $\left\{ \left\| u_{n}\right\| \right\} $ is bounded.
Now, since the embedding $H_{V,\mathrm{r}}^{1}\hookrightarrow
L_{K}^{q_{1}}+L_{K}^{q_{2}}$ is compact (see Lemma \ref{LEM: comp embedd})
and the functional $u\mapsto \int_{\mathbb{R}^{N}}K\left( \left| x\right|
\right) F\left( u\right) dx$ is of class $C^{1}$ on $%
L_{K}^{q_{1}}+L_{K}^{q_{2}}$ (see the proof of Lemma \ref{LEM: C^1}), it is
a standard exercise to conclude that $\left\{ u_{n}\right\} $ has a strongly
convergent subsequence in $H_{V,\mathrm{r}}^{1}$.%
\endproof%
\bigskip

We can now conclude the proof of Theorem \ref{THM: super general}\textbf{.}%
\bigskip

\noindent \textbf{Proof of Theorem \ref{THM: super general}. }%
Assume all the hypotheses of the theorem and assume also that $f\left( t\right) =0$
for all $t<0$. This additional hypothesis is not restrictive, since the
theorem concerns \emph{nonnegative} solutions and all its assumptions still
hold true if we replace $f\left( t\right) $ with $f\left( t\right) \chi _{%
\mathbb{R}_{+}}\left( t\right) $ (where $\chi _{\mathbb{R}_{+}}$ is the
characteristic function of $\mathbb{R}_{+}=\left( 0,+\infty \right) $).

Thanks to Lemma \ref{LEM: C^1}, the theorem is proved if we find a
nontrivial nonnegative critical point of $I:H_{V,\mathrm{r}}^{1}\rightarrow 
\mathbb{R}$.

To this end, we want to apply the Mountain-Pass Theorem \cite{Ambr-Rab}.
From (\ref{I(u)>}) of Lemma \ref{LEM: I(u)>} we deduce that, since $%
q_{1},q_{2}>2$, there exists $\rho >0$ such that 
\[
\inf_{u\in H_{V,\mathrm{r}}^{1},\,\left\| u\right\| =\rho }I\left( u\right)
>0=I\left( 0\right) . 
\]
Therefore, taking into account Lemmas \ref{LEM: C^1} and \ref{LEM:PS}, the
Mountain Pass Theorem applies if we show that $\exists \bar{u}\in H_{V,%
\mathrm{r}}^{1}$ such that $\left\| \bar{u}\right\| >\rho $ and $I\left( 
\bar{u}\right) <0$. In order to prove this, from condition $\left( \mathbf{f}%
_{6}\right) $ (which holds in any case, since it also follows from $\left( 
\mathbf{f}_{4}\right) $ and $\left( \mathbf{f}_{5}\right) $), we infer that 
\[
F\left( t\right) \geq \frac{F\left( t_{0}\right) }{t_{0}^{\theta }}t^{\theta
}\textrm{\quad for all }t\geq t_{0}. 
\]
Then we fix a radial nonnegative function $u_{0}\in C_{\mathrm{c}}^{\infty }(%
\mathbb{R}^{N}\setminus \left\{ 0\right\} )$ such that the set $\{x\in \mathbb{R}%
^{N}:u_{0}\left( x\right) \geq t_{0}\}$ has positive Lebesgue measure. We
now distinguish the case of assumptions $\left( \mathbf{f}_{4}\right) $ and $%
\left( \mathbf{f}_{5}\right) $ from the case of $K\left( \left| \cdot
\right| \right) \in L^{1}(\mathbb{R}^{N})$. In the first one, we have $F\left(
t\right) \geq 0$ for all $t\in \mathbb{R}$ (recall that $f\left( t\right) =0$
for $t<0$) and $F\left( t_{0}\right) >0$, so that for every $\lambda >1$ we
get 
\begin{eqnarray*}
\int_{\mathbb{R}^{N}}K\left( \left| x\right| \right) F\left( \lambda
u_{0}\right) dx &\geq &\int_{\left\{ \lambda u_{0}\geq t_{0}\right\}
}K\left( \left| x\right| \right) F\left( \lambda u_{0}\right) dx\geq \lambda
^{\theta }\frac{F\left( t_{0}\right) }{t_{0}^{\theta }}\int_{\left\{ \lambda
u_{0}\geq t_{0}\right\} }K\left( \left| x\right| \right) u_{0}^{\theta }dx \\
&\geq &\lambda ^{\theta }\frac{F\left( t_{0}\right) }{t_{0}^{\theta }}%
\int_{\left\{ u_{0}\geq t_{0}\right\} }K\left( \left| x\right| \right)
u_{0}^{\theta }dx\geq \lambda ^{\theta }F\left( t_{0}\right) \int_{\left\{
u_{0}\geq t_{0}\right\} }K\left( \left| x\right| \right) dx>0.
\end{eqnarray*}
Since $\theta >2$, this gives 
\begin{equation}
\lim_{\lambda \rightarrow +\infty }I\left( \lambda u_{0}\right) \leq
\lim_{\lambda \rightarrow +\infty }\left( \frac{\lambda ^{2}}{2}\left\|
u_{0}\right\| ^{2}-\lambda ^{\theta }F\left( t_{0}\right) \int_{\left\{
u_{0}\geq t_{0}\right\} }K\left( \left| x\right| \right) dx\right) =-\infty .
\label{I->-oo}
\end{equation}
If $K\left( \left| \cdot \right| \right) \in L^{1}(\mathbb{R}^{N})$, assumption 
$\left( \mathbf{f}_{6}\right) $ still gives $F\left( t_{0}\right) >0$ and
from (\ref{|F| <}) we infer that 
\[
F\left( t\right) \geq -\tilde{M}\min \left\{
t_{0}^{q_{1}},t_{0}^{q_{2}}\right\} \qquad \textrm{for all }0\leq t\leq t_{0}. 
\]
Therefore, arguing as above about the integral over $\left\{ \lambda
u_{0}\geq t_{0}\right\} $, for every $\lambda >1$ we obtain 
\begin{eqnarray*}
\int_{\mathbb{R}^{N}}K\left( \left| x\right| \right) F\left( \lambda
u_{0}\right) dx &=&\int_{\left\{ \lambda u_{0}<t_{0}\right\} }K\left( \left|
x\right| \right) F\left( \lambda u_{0}\right) dx+\int_{\left\{ \lambda
u_{0}\geq t_{0}\right\} }K\left( \left| x\right| \right) F\left( \lambda
u_{0}\right) dx \\
&\geq &-\tilde{M}\min \left\{ t_{0}^{q_{1}},t_{0}^{q_{2}}\right\}
\int_{\left\{ \lambda u_{0}<t_{0}\right\} }K\left( \left| x\right| \right)
dx+\lambda ^{\theta }F\left( t_{0}\right) \int_{\left\{ u_{0}\geq
t_{0}\right\} }K\left( \left| x\right| \right) dx,
\end{eqnarray*}
which implies 
\[
\lim_{\lambda \rightarrow +\infty }I\left( \lambda u_{0}\right) \leq
\lim_{\lambda \rightarrow +\infty }\left( \frac{\lambda ^{2}}{2}\left\|
u_{0}\right\| ^{2}+\tilde{M}\min \left\{ t_{0}^{q_{1}},t_{0}^{q_{2}}\right\}
\left\| K\right\| _{L^{1}(\mathbb{R}^{N})}-\lambda ^{\theta }F\left(
t_{0}\right) \int_{\left\{ u_{0}\geq t_{0}\right\} }K\left( \left| x\right|
\right) dx\right) =-\infty . 
\]
So, in any case, we can take $\bar{u}=\lambda u_{0}$ with $\lambda $
sufficiently large and the Mountain-Pass Theorem provides the existence of a
nontrivial critical point $u\in H_{V,\mathrm{r}}^{1}$ for $I$. Since $%
f\left( t\right) =0$ for $t<0$ implies $I^{\prime }\left( u\right)
u_{-}=-\left\| u_{-}\right\| ^{2}$ (where $u_{-}\in H_{V,\mathrm{r}}^{1}$ is
the negative part of $u$), one concludes that $u_{-}=0$, i.e., $u$ is
nonnegative.%
\endproof%
\bigskip

For concluding also the proof of Theorem \ref{THM: sub general}, we prove
one more lemma\textbf{.}

\begin{lemma}
\label{LEM:inf<0}If $\left( \mathbf{f}_{8}\right) $ holds, then the
functional $I:H_{V,\mathrm{r}}^{1}\rightarrow \mathbb{R}$ takes negative values.
\end{lemma}

\proof%
By assumption $\left( \mathbf{f}_{8}\right) $, fix $m>0$ and $t_{0}>0$ such
that $F\left( t\right) \geq mt^{\theta }$ for all $0\leq t\leq t_{0}$. Fix a
nonzero radial function $u_{0}\in C_{\mathrm{c}}^{\infty }(\mathbb{R}%
^{N}\setminus \left\{ 0\right\} )$ such that $0\leq u_{0}\leq t_{0}$. Then,
for every $0<\lambda <1$ we get that $\lambda u_{0}\in H_{V,\mathrm{r}}^{1}$
satisfies $0\leq \lambda u_{0}\leq t_{0}$ and therefore 
\[
I\left( \lambda u_{0}\right) =\frac{1}{2}\left\| \lambda u_{0}\right\|
^{2}-\int_{\mathbb{R}^{N}}K\left( \left| x\right| \right) F\left( \lambda
u_{0}\right) dx\leq \frac{\lambda ^{2}}{2}\left\| u_{0}\right\| ^{2}-\lambda
^{\theta }m\int_{\mathbb{R}^{N}}K\left( \left| x\right| \right) u_{0}^{\theta
}dx, 
\]
where $\int_{\mathbb{R}^{N}}K\left( \left| x\right| \right) u_{0}^{\theta }dx>0$
(recall that $K>0$ everywhere) .Since $\theta <2$, this implies $I\left(
\lambda u_{0}\right) <0$ for $\lambda $ sufficiently small.%
\endproof%
\bigskip

\noindent \textbf{Proof of Theorem \ref{THM: sub general}. }%
Assume all the hypotheses of the theorem and assume also that $f$ is odd. This additional
hypothesis is not restrictive, since the theorem concerns \emph{nonnegative}
solutions and all its assumptions still hold true if we replace $f\left(
t\right) $ with $f\left( \left| t\right| \right) \sgn\left( t\right) $
(where $\sgn$ is the sign function).

Since $q_{1},q_{2}\in \left( 1,2\right) $, the inequality (\ref{I(u)>}) of
Lemma \ref{LEM: I(u)>} readily implies that the functional $I:H_{V,\mathrm{r}%
}^{1}\rightarrow \mathbb{R}$ is bounded from below and coercive, so that 
\[
\mu :=\inf_{u\in H_{V,\mathrm{r}}^{1}}I\left( u\right) 
\]
is a finite value. Therefore, thanks to Lemma \ref{LEM: C^1}, the theorem is
proved if we show that $\mu $ is attained by a nonnegative minimizer, which
cannot be trivial, since $I\left( 0\right) =0$ and $\mu <0$, by Lemma \ref
{LEM:inf<0}.

To this end, let $\left\{ u_{n}\right\} $ be any minimizing sequence for $%
\mu $. Since $f$ is odd, $I\left( u\right) $ is even and therefore $\left\{
\left| u_{n}\right| \right\} $ is still a minimizing sequence, so that, up
to replacing $u_{n}$ with $\left| u_{n}\right| $, we may assume $u_{n}\geq 0$%
. Since $\left\{ u_{n}\right\} $ is bounded in $H_{V,\mathrm{r}}^{1}$ (by
the coercivity of $I$) and the embedding $H_{V,\mathrm{r}}^{1}%
\hookrightarrow L_{K}^{q_{1}}+L_{K}^{q_{2}}$ is compact (by Lemma \ref{LEM:
comp embedd}), up to a subsequence we can assume that there exists $u\in
H_{V,\mathrm{r}}^{1}$ such that: 
\[
u_{n}\rightharpoonup u\quad \textrm{in }H_{V,\mathrm{r}}^{1},\quad
u_{n}\rightarrow u\quad \textrm{in }L_{K}^{q_{1}}+L_{K}^{q_{2}},\quad
u_{n}\rightarrow u\quad \textrm{almost everywhere in }\mathbb{R}^{N} 
\]
(the almost everywhere convergence follows, for instance, from the
continuous embedding $H_{V,\mathrm{r}}^{1}\hookrightarrow D^{1,2}(\mathbb{R}%
^{N})$ and the fact that, up to a subsequence, weak convergence in $D^{1,2}(%
\mathbb{R}^{N})$ implies almost everywhere convergence). Then $u_{n}\geq 0$
implies $u\geq 0$ and, thanks to the weak lower semi-continuity of the norm
and to the continuity of the functional $v\mapsto \int_{\mathbb{R}^{N}}K\left(
\left| x\right| \right) F\left( v\right) dx$ on $L_{K}^{q_{1}}+L_{K}^{q_{2}}$
(see the proof of Lemma \ref{LEM: C^1} above), $u$ satisfies 
\[
\left\| u\right\| ^{2}\leq \liminf_{n\rightarrow \infty }\left\|
u_{n}\right\| ^{2}\quad \textrm{and}\quad \int_{\mathbb{R}^{N}}K\left( \left|
x\right| \right) F\left( u\right) dx=\lim_{n\rightarrow \infty }\int_{\mathbb{R}%
^{N}}K\left( \left| x\right| \right) F\left( u_{n}\right) dx. 
\]
This implies 
\[
I\left( u\right) =\frac{1}{2}\left\| u\right\| ^{2}-\int_{\mathbb{R}%
^{N}}K\left( \left| x\right| \right) F\left( u\right) dx\leq
\lim_{n\rightarrow \infty }\left( \frac{1}{2}\left\| u_{n}\right\|
^{2}-\int_{\mathbb{R}^{N}}K\left( \left| x\right| \right) F\left( u_{n}\right)
dx\right) =\mu 
\]
and therefore we conclude $I\left( u\right) =\mu $.%
\endproof%

\section{Existence of a ground state \label{SEC: GS}}

In this section we give a version of Theorem \ref{THM: super general} which
ensures the existence of a \emph{radial ground state} of Eq. (\ref{Eq.}), by
which, assuming that the Euler functional $I$ defined in (\ref{I:=}) is of
class $C^{1}$ on $H_{V,\mathrm{r}}^{1}$ (as in Lemma \ref{LEM: C^1}), we
mean a radial solution $u\neq 0$ such that 
\[
I\left( u\right) =\min_{v\in \mathcal{N}}I\left( v\right) \textrm{\quad
where\quad }\mathcal{N}:=\left\{ v\in H_{V,\mathrm{r}}^{1}\setminus \left\{
0\right\} :I^{\prime }\left( v\right) v=0\right\} 
\]
($\mathcal{N}$ is the Nehari manifold). As $I$ is often called the
``action'' or ``energy'' functional associated to the equation, a radial
ground state $u$ is in fact a \emph{least action} or \emph{least energy}
solution (among the nontrivial radial ones), since every radial solution $%
v\neq 0$ belongs to $\mathcal{N}$ and therefore $I\left( u\right) \leq
I\left( v\right) $. Of course, the solution found in Theorem \ref{THM: sub
general} is itself a radial ground state, since it is a global minimizer of $%
I$ on $H_{V,\mathrm{r}}^{1}$.

The result we will prove is the following theorem. Observe that all its
assumptions are satisfied by both the first two nonlinearities of Example 
\ref{EX:f}, with $q_{1},q_{2}>2$, $q_{1}\in \mathcal{I}_{1}$, $q_{2}\in \mathcal{I}_{2}$.

\begin{theorem}
\label{THM: GS}Under the same assumptions of the first part of Theorem \ref
{THM: super general} (i.e., the part with $\left( \mathbf{f}_{3}\right) $-$%
\left( \mathbf{f}_{5}\right) $), if $f$ also satisfies

\begin{itemize}
\item[$\left( \mathbf{f}_{9}\right) $]  
$\displaystyle\,\frac{f\left( t\right) }{t}$ is a strictly increasing function on $\left(
0,+\infty \right) ,$
\end{itemize}

\noindent then Eq. (\ref{Eq.}) has a nonnegative radial ground state.
\end{theorem}

In proving Theorem \ref{THM: GS}, we will use the following variant of the
Mountain Pass Theorem and an adaptation of well known arguments involving
the Nehari manifold (see e.g. \cite[Chapter 4]{Willem}).

\begin{lemma}
\label{LEM:minimax}Let $X$ be a real Banach space and let $J\in C^{1}(X;\mathbb{%
R})$. Assume that there exist $\rho >0$ and $\overline{v}\in X$ such that 
\[
\inf_{v\in X,\,\left\| v\right\| _{X}=\rho }J\left( v\right) >\inf_{v\in
X,\,\left\| v\right\| _{X}\leq \rho }J\left( v\right) =J\left( 0\right)
=0>J\left( \overline{v}\right) \quad \textrm{and}\quad \left\| \overline{v}%
\right\| _{X}>\rho .
\]
If $J$ satisfies the Palais-Smale condition, then the minimax level 
\begin{equation}
c:=\inf_{\gamma \in \Gamma }\max_{0\leq t\leq 1}J\left( \gamma \left(
t\right) \right) >0,\quad \Gamma :=\left\{ \gamma \in C\left( \left[
0,1\right] ;X\right) :\gamma \left( 0\right) =0,\,J\left( \gamma \left(
1\right) \right) <0\right\} ,  \label{minimax}
\end{equation}
is a critical value for $J$.
\end{lemma}

\proof%
It follows for instance from \cite[Theorem 2.9]{Willem} (apply with $%
M=\left[ 0,1\right] $, $M_{0}=\left\{ 0,1\right\} $ and $\Gamma _{0}=\left\{
\gamma \in C\left( M_{0};X\right) :\gamma \left( 0\right) =0,\,J\left(
\gamma \left( 1\right) \right) <0\right\} $).%
\endproof%
\bigskip

\noindent \textbf{Proof of Theorem \ref{THM: GS}. }%
Assume all the hypotheses of the theorem. As in the proof of Theorem \ref{THM: super
general}, we assume additionally that $f\left( t\right) =0$ for all $t<0$
and we deduce that $I$ satisfies the assumptions of Lemma \ref{LEM:minimax}
(with $X=H_{V,\mathrm{r}}^{1}$ and $J=I$), so that there exists $u\in H_{V,%
\mathrm{r}}^{1}$ such that $I^{\prime }\left( u\right) =0$, $u\geq 0$ and $%
I\left( u\right) =c>0$, where $c$ is the minimax level (\ref{minimax}).
Since $u\neq 0$ is a critical point for $I$, we have that $u\in \mathcal{N}$
and therefore $c=I\left( u\right) \geq \nu :=\inf_{v\in \mathcal{N}}I\left(
v\right) $. Hence the theorem is proved if we show that $\nu \geq c$, which
implies $I\left( u\right) =\nu $ with $u\in \mathcal{N}$.

To this end, we take any $v\in \mathcal{N}$ and first observe that it cannot
be $v\leq 0$ almost everywhere. Otherwise, since $v\in \mathcal{N}$ implies $%
\left\| v\right\| ^{2}=\int_{\mathbb{R}^{N}}K\left( \left| x\right| \right)
f\left( v\right) v\,dx$ (recall Lemma \ref{LEM: C^1}) and we have $f\left(
t\right) =0$ for $t\leq 0$, we would get the contradiction $\left\|
v\right\| =0$.

Then we show that $I\left( v\right) =\max_{t\geq 0}I\left( tv\right) >0$.
For this, we define $g_{v}\left( t\right) :=I\left( tv\right) $, $t\geq 0$,
and argue by the following three steps.

\begin{itemize}
\item  $t=1$ is a critical point for $g_{v}$, since $g_{v}^{\prime }\left(
t\right) =I^{\prime }\left( tv\right) v$ for all $t$ and $I^{\prime }\left(
v\right) v=0$.

\item  $t=1$ is the only critical point of $g_{v}$ on $\left( 0,+\infty
\right) $. Indeed, if $t_{2}>t_{1}>0$ are critical points for $g_{v}$, then
we have $I^{\prime }\left( t_{1}v\right) v=I^{\prime }\left( t_{2}v\right)
v=0$, i.e., 
\[
t_{1}\left\| v\right\| ^{2}-\int_{\mathbb{R}^{N}}K\left( \left| x\right|
\right) f\left( t_{1}v\right) v\,dx=t_{2}\left\| v\right\| ^{2}-\int_{\mathbb{R}%
^{N}}K\left( \left| x\right| \right) f\left( t_{2}v\right) v\,dx=0
\]
(recall Lemma \ref{LEM: C^1}), which implies 
\[
0=\int_{\mathbb{R}^{N}}K\left( \left| x\right| \right) \left( \frac{f\left(
t_{2}v\right) }{t_{2}}-\frac{f\left( t_{1}v\right) }{t_{1}}\right)
v\,dx=\int_{\left\{ v>0\right\} }K\left( \left| x\right| \right) \left( 
\frac{f\left( t_{2}v\right) }{t_{2}v}-\frac{f\left( t_{1}v\right) }{t_{1}v}%
\right) v^{2}\,dx
\]
(recall that $f\left( t\right) =0$ for $t\leq 0$). The last integrand is
nonnegative by assumption $\left( \mathbf{f}_{9}\right) $ and therefore we
get $f\left( t_{2}v\right) /\left( t_{2}v\right) =f\left( t_{1}v\right)
/\left( t_{1}v\right) $ almost everywhere on $\left\{ x:v\left( x\right)
>0\right\} $ (because $K\left( \left| x\right| \right) ,v^{2}\left( x\right)
>0$ and the set $\left\{ x:v\left( x\right) >0\right\} $ has positive
measure). Again by assumption $\left( \mathbf{f}_{9}\right) $, this implies $%
t_{2}v=t_{1}v$ almost everywhere on $\left\{ x:v\left( x\right) >0\right\} $
and therefore we conclude $t_{2}=t_{1}$.

\item  $g_{v}$ has a maximum point on $\left( 0,+\infty \right) $, in which $%
g_{v}>0$. Indeed, letting $\delta >0$ be such that the set $\left\{
x:v\left( x\right) \geq \delta \right\} $ has positive measure (which exists
because $\left\{ x:v\left( x\right) >0\right\} $ has positive measure), we
have $F\left( \delta \right) >0$ (because $\left( \mathbf{f}_{4}\right) $, $%
\left( \mathbf{f}_{9}\right) $ and $f\left( 0\right) =0$ imply $F\left(
t\right) >0$ for $t>0$) and from $\left( \mathbf{f}_{4}\right) $ we deduce
that $F\left( t\right) \geq F\left( \delta \right) \delta ^{-\theta
}t^{\theta }$ for all $t\geq \delta $, so that, arguing as for (\ref{I->-oo}%
), we get 
\begin{equation}
\lim_{t\rightarrow +\infty }I\left( tv\right) \leq \lim_{t\rightarrow
+\infty }\left( \frac{t^{2}}{2}\left\| v\right\| ^{2}-t^{\theta }F\left(
\delta \right) \int_{\left\{ v\left( x\right) \geq \delta \right\} }K\left(
\left| x\right| \right) dx\right) =-\infty .  \label{I->-oo bis}
\end{equation}
This, together with $g_{v}\left( 0\right) =0$ and $g_{v}\left( t\right) >0$
for $t>0$ small enough (which follows from (\ref{I(u)>}) of Lemma \ref{LEM:
I(u)>}, where $q_{1},q_{2}>2$), yields the claim.
\end{itemize}

\noindent As a result, the maximum point of $g_{v}$ on $\left( 0,+\infty
\right) $ must be the unique critical point $t=1$ and we conclude that $%
I\left( v\right) =\max_{t\geq 0}I\left( tv\right) >0$.

Now, using (\ref{I->-oo bis}) again, we observe that there exists $t_{v}>0$
satisfying $I\left( tv\right) <0$ for all $t\geq t_{v}$, so that the path $%
\gamma _{v}\left( t\right) :=tt_{v}v$, $t\in \left[ 0,1\right] $, is such
that $\gamma _{v}\in \Gamma $ and 
\[
\max_{0\leq t\leq 1}I\left( \gamma _{v}\left( t\right) \right) =\max_{0\leq
t\leq 1}I\left( tt_{v}v\right) =\max_{0\leq t\leq t_{v}}I\left( tv\right)
=\max_{t\geq 0}I\left( tv\right) =I\left( v\right) . 
\]
Hence 
\[
I\left( v\right) =\max_{0\leq t\leq 1}I\left( \gamma _{v}\left( t\right)
\right) \geq \inf_{\gamma \in \Gamma }\max_{0\leq t\leq 1}I\left( \gamma
\left( t\right) \right) =c 
\]
and therefore, as $v\in \mathcal{N}$ is arbitrary, we get $\nu \geq c$.%
\endproof%

}


\begin{thebibliography}{99}
\bibitem{Alves-Souto-13}  \textsc{Alves C.O.}, \textsc{Souto M.A.S.}, \emph{%
Existence of solutions for a class of nonlinear Schr\"{o}dinger equations
with potential vanishing at infinity}, J. Differential Equations \textbf{254}
(2013), 1977-1991.

\bibitem{Ambr-Fel-Malch}  \textsc{Ambrosetti A.}, \textsc{Felli V.}, \textsc{%
Malchiodi A.}, \emph{Ground states of nonlinear Schr\"{o}dinger equations
with potentials vanishing at infinity}, J. Eur. Math. Soc. \textbf{7}
(2005), 117-144.

\bibitem{Ambr-Rab}  \textsc{Ambrosetti A.}, \textsc{Rabinowitz P.H.}, \emph{%
Dual variational methods in critical point theory and applications}, J.
Funct. Anal. \textbf{14} (1973), 349-381.

\bibitem{Azz-D'Av-Pomp14}  \textsc{Azzollini A.}, \textsc{d'Avenia P.}, 
\textsc{Pomponio A.}, \emph{Quasilinear elliptic equations in }$\mathbb{R}^{N}$%
\emph{\ via variational methods and Orlicz-Sobolev embeddings}, Calc. Var.
Partial Differential Equations \textbf{49} (2014), 197-213.

\bibitem{BBR1}  \textsc{Badiale M.}, \textsc{Benci V.}, \textsc{Rolando S.}, 
\emph{Solitary waves: physical aspects and mathematical results}, Rend. Sem.
Math. Univ. Pol. Torino \textbf{62} (2004), 107-154.

\bibitem{BGRnonex}  \textsc{Badiale M.}, \textsc{Guida M.}, \textsc{Rolando
S.}, \emph{A nonexistence result for a nonlinear elliptic equation with
singular and decaying potential}, Commun. Contemp. Math., in press (DOI:
10.1142/S0219199714500242).

\bibitem{BGR-p1}  \textsc{Badiale M.}, \textsc{Guida M.}, \textsc{Rolando S.}%
, \emph{Compactness and existence results in weighted Sobolev spaces of
radial functions. Part I: Compactness}, Cal. Var. Partial Differential
Equations, in press (DOI: 10.1007/s00526-015-0817-2).

\bibitem{BGR-p2}  \textsc{Badiale M.}, \textsc{Guida M.}, \textsc{Rolando S.}%
, \emph{Compactness and existence results in weighted Sobolev spaces of
radial functions. Part II: Existence}, in preparation.

\bibitem{BPR}  \textsc{Badiale M.}, \textsc{Pisani L.}, \textsc{Rolando S.}, 
\emph{Sum of weighted Lebesgue spaces and nonlinear elliptic equations},
NoDEA, Nonlinear Differ. Equ. Appl. \textbf{18} (2011), 369-405.

\bibitem{BR}  \textsc{Badiale M.}, \textsc{Rolando S.}, \emph{Elliptic
problems with singular potentials and double-power nonlinearity}, Mediterr.
J. Math. \textbf{2} (2005), 417-436.

\bibitem{BR TMA}  \textsc{Badiale M., Rolando S.}, \emph{Nonlinear elliptic
equations with subhomogeneous potentials}, Nonlinear Anal. \textbf{72}
(2010), 602-617.

\bibitem{BFmonograph}  \textsc{Benci V.}, \textsc{Fortunato D.}, \emph{%
Variational methods in nonlinear field equations. Solitary waves,
hylomorphic solitons and vortices}, Springer Monographs in Mathematics,
Springer, Cham, 2014.

\bibitem{Be-Gr-Mic}  \textsc{Benci V.}, \textsc{Grisanti C.R.}, \textsc{%
Micheletti A.M.}, \emph{Existence and non existence of the ground state
solution for the nonlinear Schr\"{o}dinger equation with }$V\left( \infty
\right) =0$, Topol. Methods Nonlinear Anal. \textbf{26} (2005), 203-220.

\bibitem{Be-Gr-Mic 2}  \textsc{Benci V.}, \textsc{Grisanti C.R.}, \textsc{%
Micheletti A.M.}, \emph{Existence of solutions for the nonlinear Schr\"{o}%
dinger equation with }$V\left( \infty \right) =0$, Contributions to
nonlinear analysis, Progr. Nonlinear Differential Equations Appl., vol. 66,
Birkh\"{a}user, Basel, 2006.

\bibitem{Beres-Lions}  \textsc{Berestycki H.}, \textsc{Lions P.L.}, \emph{%
Nonlinear Scalar Field Equations, I - II}, Arch. Rational Mech. Anal. 
\textbf{82} (1983), 313-379.

\bibitem{Catrina nonex}  \textsc{Catrina F.}, \emph{Nonexistence of positive
radial solutions for a problem with singular potential}, Adv. Nonlinear
Anal. \textbf{3} (2014), 1-13.

\bibitem{Floer-Wein}  \textsc{Floer A.}, \textsc{Weinstein A.}, \emph{%
Nonspreading wave packets for the cubic Schr\"{o}dinger equation with a
bounded potential}, J.\ Funct. Anal. \textbf{69} (1986), 397-408.

\bibitem{GR-boundedPS}  \textsc{Guida M.}, \textsc{Rolando S.}, \emph{On the
existence of bounded Palais-Smale sequences and applications to nonlinear
equations without superlinearity assumptions}, in preparation.

\bibitem{Rabi92}  \textsc{Rabinowitz P.H.}, \emph{On a class of nonlinear
Schr\"{o}dinger equations}, Z. Angew. Math. Phys. \textbf{43} (1992),
270-291.

\bibitem{Strauss}  \textsc{Strauss W.A.}, \emph{Existence of solitary waves
in higher dimensions}, Comm. Math. Phys. \textbf{55} (1977), 149-172.

\bibitem{SuTian12}  \textsc{Su J., Tian R.}, \emph{Weighted Sobolev type
embeddings and coercive quasilinear elliptic equations on }$\mathbb{R}^{N}$,
Proc. Amer. Math. Soc. \textbf{140} (2012), 891-903.

\bibitem{Su-Wang-Will p}  \textsc{Su J., Wang Z.-Q., Willem M.}, \emph{%
Weighted Sobolev embedding with unbounded and decaying radial potentials},
J. Differential Equations \textbf{238} (2007), 201-219.

\bibitem{Willem}  \textsc{Willem M.}, \emph{Minimax Theorems}, Progress in
Nonlinear Differential Equations and their Applications, vol. 24, Birkh\"{a}%
user, Boston, 1996.

\bibitem{YangY}  \textsc{Yang Y.}, \emph{Solitons in field theory and
nonlinear analysis}, Springer Monographs in Mathematics, Springer-Verlag,
New York, 2001.

\bibitem{Zhang-Fu15}  \textsc{Zhang G.}, \textsc{Fu H.}, \emph{Ground states
for a modified capillary surface equation in weighted Orlicz-Sobolev space},
Electron. J. Differential Equations \textbf{2015} (2015), No. 85, 1-18.
\bigskip \bigskip \pagebreak 

\end{thebibliography}
\end{document}